\documentclass[11pt,a4paper,leqno,noamsfonts]{amsart}
\usepackage[english]{babel}
\usepackage[dvipsnames]{xcolor}
\usepackage{graphicx,pifont}
\usepackage[utopia]{mathdesign}
\usepackage[a4paper,top=3.2cm,bottom=3.2cm,left=3.5cm,right=3.5cm,marginparwidth=60pt]{geometry}
\usepackage[utf8]{inputenc}
\usepackage{braket,caption,comment,mathtools,stmaryrd,enumitem} 
\usepackage[usestackEOL]{stackengine}
\usepackage{multirow,booktabs,microtype,relsize,soul}
\usepackage[foot]{amsaddr} 
\usepackage[colorlinks,bookmarks,backref=page]{hyperref} %
\hypersetup{colorlinks,%
        citecolor=britishracinggreen,%
        filecolor=black,%
        linkcolor=cobalt,%
        urlcolor=purple}
\numberwithin{equation}{section}

\usepackage[capitalise]{cleveref}

\allowdisplaybreaks

\Crefname{conjecture}{Conjecture}{Conjectures}
\Crefname{prop}{Proposition}{Propositions}

\newcounter{x}
\newcounter{y}
\newcounter{z}

\newcommand\xaxis{210}
\newcommand\yaxis{-30}
\newcommand\zaxis{90}

\newcommand\topside[3]{
  \fill[fill=yellow, draw=black,shift={(\xaxis:#1)},shift={(\yaxis:#2)},
  shift={(\zaxis:#3)}] (0,0) -- (30:1) -- (0,1) --(150:1)--(0,0);
}
\newcommand\topsideone[3]{
  \fill[fill=magenta, draw=black,shift={(\xaxis:#1)},shift={(\yaxis:#2)},
  shift={(\zaxis:#3)}] (0,0) -- (30:1) -- (0,1) --(150:1)--(0,0);
}

\newcommand\leftside[3]{
  \fill[fill=red, draw=black,shift={(\xaxis:#1)},shift={(\yaxis:#2)},
  shift={(\zaxis:#3)}] (0,0) -- (0,-1) -- (210:1) --(150:1)--(0,0);
}
\newcommand\leftsideone[3]{
  \fill[fill=purple, draw=black,shift={(\xaxis:#1)},shift={(\yaxis:#2)},
  shift={(\zaxis:#3)}] (0,0) -- (0,-1) -- (210:1) --(150:1)--(0,0);
}

\newcommand\rightside[3]{
  \fill[fill=blue, draw=black,shift={(\xaxis:#1)},shift={(\yaxis:#2)},
  shift={(\zaxis:#3)}] (0,0) -- (30:1) -- (-30:1) --(0,-1)--(0,0);
}
\newcommand\rightsideone[3]{
  \fill[fill=brown, draw=black,shift={(\xaxis:#1)},shift={(\yaxis:#2)},
  shift={(\zaxis:#3)}] (0,0) -- (30:1) -- (-30:1) --(0,-1)--(0,0);
}

\newcommand\cube[3]{
  \topside{#1}{#2}{#3} \leftside{#1}{#2}{#3} \rightside{#1}{#2}{#3}
}

\newcommand\cubeone[3]{
  \topsideone{#1}{#2}{#3} \leftsideone{#1}{#2}{#3} \rightsideone{#1}{#2}{#3}
}

\newcommand\planepartition[1]{
 \setcounter{x}{-1}
  \foreach \a in {#1} {
    \addtocounter{x}{1}
    \setcounter{y}{-1}
    \foreach \b in \a {
      \addtocounter{y}{1}
      \setcounter{z}{-1}
      \foreach \c in {1,...,\b} {
        \addtocounter{z}{1}
        \cube{\value{x}}{\value{y}}{\value{z}}
      }
    }
  }
}

\newcommand\planepartitionone[1]{
 \setcounter{x}{-1}
  \foreach \a in {#1} {
    \addtocounter{x}{1}
    \setcounter{y}{-1}
    \foreach \b in \a {
      \addtocounter{y}{1}
      \setcounter{z}{-1}
      \foreach \c in {1,...,\b} {
        \addtocounter{z}{1}
        \cube{\value{x}}{\value{y}}{\value{z}}
      }
    }
  }
      \setcounter{x}{0}
      \setcounter{y}{0}
      \setcounter{z}{-1}
      \cubeone{\value{x}}{\value{y}}{\value{z}}
      \setcounter{x}{1}
      \setcounter{y}{-1}
      \setcounter{z}{-1}
      \cubeone{\value{x}}{\value{y}}{\value{z}}
}

\newcommand\planepartitiondue[1]{
 \setcounter{x}{-1}
  \foreach \a in {#1} {
    \addtocounter{x}{1}
    \setcounter{y}{-1}
    \foreach \b in \a {
      \addtocounter{y}{1}
      \setcounter{z}{-1}
      \foreach \c in {1,...,\b} {
        \addtocounter{z}{1}
        \cube{\value{x}}{\value{y}}{\value{z}}
      }
    }
  }
      \setcounter{x}{0}
      \setcounter{y}{-1}
      \setcounter{z}{-1}
      \cubeone{\value{x}}{\value{y}}{\value{z}}
      \setcounter{x}{0}
      \setcounter{y}{0}
      \setcounter{z}{-1}
      \cubeone{\value{x}}{\value{y}}{\value{z}}
      \setcounter{x}{1}
      \setcounter{y}{-1}
      \setcounter{z}{-1}
      \cubeone{\value{x}}{\value{y}}{\value{z}}
}
\DeclareSymbolFont{cmarrows}{OMS}{cmsy}{m}{n}
\SetSymbolFont{cmarrows}{bold}{OMS}{cmsy}{b}{n}
\DeclareMathSymbol{\cmminus}{\mathbin}{cmarrows}{"00}

\DeclareMathSymbol{\leftrightarrow}{\mathrel}{cmarrows}{"24}
\DeclareMathSymbol{\leftarrow}{\mathrel}{cmarrows}{"20}
   
\DeclareMathSymbol{\rightarrow}{\mathrel}{cmarrows}{"21}
   \let\to=\rightarrow
\DeclareMathSymbol{\mapstochar}{\mathrel}{cmarrows}{"37}
   \def\mapsto{\mapstochar\rightarrow}
   
\DeclareSymbolFont{usualmathcal}{OMS}{cmsy}{m}{n}
\DeclareSymbolFontAlphabet{\mathcal}{usualmathcal}
\DeclareMathAlphabet\BCal{OMS}{cmsy}{b}{n}

\makeatletter
\newcommand{\mylabel}[2]{#2\def\@currentlabel{#2}\label{#1}}
\makeatother

\definecolor{cornellred}{rgb}{0.7, 0.11, 0.11}
\definecolor{britishracinggreen}{rgb}{0.0, 0.26, 0.15}
\definecolor{cobalt}{rgb}{0.0, 0.28, 0.67}
 
\usepackage{amsmath,float,soul}

\DeclareMathOperator{\IB}{\mathbf{IB}}
\DeclareMathOperator{\Bor}{\mathbf{B}}

\newcommand{\Coef}{\mathsf{Coef}}

\newcommand{\onto}{\twoheadrightarrow}

\newcommand{\boldit}[1]{\boldsymbol{#1}}

\newcommand{\CCoh}{\mathscr{C}\kern-0.25em {o}\kern-0.2em{h}}

\DeclareMathOperator{\ba}{{\boldit{a}}}
\DeclareMathOperator{\by}{{\boldit{y}}}

\DeclareMathOperator{\Soc}{\mathrm{Soc}}

\DeclareMathOperator{\Hilb}{Hilb}

\DeclareMathOperator{\aut}{aut}

\DeclareMathOperator{\Quot}{Quot}

\newcommand{\BA}{{\mathbb{A}}}

\newcommand{\BC}{{\mathbb{C}}}

\newcommand{\BL}{{\mathbb{L}}}

\newcommand{\BN}{{\mathbb{N}}}

\newcommand{\BQ}{{\mathbb{Q}}}
\newcommand{\BR}{{\mathbb{R}}}

\newcommand{\BZ}{{\mathbb{Z}}}

\newcommand{\FS}{{\mathfrak{S}}}

\usepackage[all]{xy}
\usepackage{tikz}
\usepackage{tikz-cd}
\usepackage{adjustbox}
\usepackage{rotating}
\usepackage{comment}

\usetikzlibrary{matrix,shapes,intersections,arrows,decorations.pathmorphing}
\tikzset{commutative diagrams/.cd,
mysymbol/.style={start anchor=center,end anchor=center,draw=none}}

\tikzset{
shift up/.style={
to path={([yshift=#1]\tikztostart.east) -- ([yshift=#1]\tikztotarget.west) \tikztonodes}
}
}

\usepackage{youngtab} 

\usepackage{ytableau} 
 \ytableausetup{boxsize=1.09em}
       \ytableausetup{boxframe=0.02em}\ytableausetup{aligntableaux=bottom}

\theoremstyle{definition}

\newtheorem*{lemma*}{Lemma}
\newtheorem*{theorem*}{Theorem}
\newtheorem*{example*}{Example}
\newtheorem*{fact*}{Fact}
\newtheorem*{notation*}{Notation}
\newtheorem*{definition*}{Definition}
\newtheorem*{prop*}{Proposition}
\newtheorem*{remark*}{Remark}
\newtheorem*{corollary*}{Corollary}

\newtheorem*{conventions*}{Conventions}

\newtheorem{definition}{Definition}[section]

\newtheorem{example}[definition]{Example}

\newtheorem{notation}[definition]{Notation}

\newtheorem{remark}[definition]{Remark}

\newtheoremstyle{thm} 
        {3mm}
        {3mm}
        {\slshape}
        {0mm}
        {\bfseries}
        {.}
        {1mm}
        {}
\theoremstyle{thm}
\newtheorem{theorem}[definition]{Theorem}
\newtheorem{corollary}[definition]{Corollary}
\newtheorem{lemma}[definition]{Lemma}
\newtheorem{prop}[definition]{Proposition}
\newtheorem{conjecture}[definition]{Conjecture}

\newtheorem{thm}{Theorem}

\newtheorem*{Acknowledgments*}{Acknowledgments}

\usepackage{todonotes}
\usepackage{marginnote}
\usepackage{xpatch}
\usepackage{tcolorbox}

\makeatletter
\patchcmd{\@todonotes@drawMarginNoteWithLine}{\marginpar}{\marginnote}{}{}
\makeatother

\makeatletter
\newenvironment{proofof}[1]{\par
  \pushQED{\qed}%
  \normalfont \topsep6\p@\@plus6\p@\relax
  \trivlist
  \item[\hskip3\labelsep
        \itshape
    Proof of #1\@addpunct{.}]\ignorespaces
}{%
  \popQED\endtrivlist\@endpefalse
}
\makeatother

\newcommand{\ancillary}{\href{https://www.rimosco.it/generatingfunctions.txt}{{\color{blue}generatingfunctions.txt}}}

\title[Enumeration of partitions via socle reduction]{Enumeration of partitions via socle reduction}
\author{Michele Graffeo, Sergej Monavari, Riccardo Moschetti, Andrea T. Ricolfi}

\begin{document}

\begin{abstract}
We study the enumeration problem of higher dimensional partitions, a natural generalisation of classical integer partitions. We show that their counting problem is equivalent to the enumeration of simpler classes of higher dimensional partitions,  satisfying suitable constraints on their embedding dimension and socle type. We provide exact formulas for the generating functions of several infinite families of such partitions, and design a procedure enumerating them in the general case. As a proof of concept, we  determine the number of partitions of size up to 30 in any dimension.
\end{abstract}

\dedicatory{In memory of Percy A.~MacMahon (1854--1929), on the occasion of his $\binom{19}{2}$-th birthday.}

\maketitle
{\hypersetup{linkcolor=black}\tableofcontents}

\section{Introduction}
\subsection{Overview}
\emph{Higher dimensional partitions} are classical combinatorial objects introduced by MacMahon over a century ago \cite{MacMahon1,MacMahon2}. For $n,d \in \BZ_{\geqslant 0}$, an $(n-1)$-dimensional partition of \emph{size} $d$ can be graphically visualised as an arrangement of $d$ boxes in $\BR^n$. Let $p_d^n$ be the number of such partitions. Based on small experimental data, MacMahon conjectured the elegant formula
\begin{equation}
\label{eqn: MacMahon conj intro}
  \sum_{d\geqslant 0}p_d^nt^d=\prod_{m>0} \left(1-t^m\right)^{-\binom{m+n-3}{n-2}},  
\end{equation}
which does hold for $n=2,3$. The case $n=2$ was known to Euler, whereas for the case $n=3$ we refer to MacMahon's original work \cite{MacMahon1,MacMahon2}.
The identity \eqref{eqn: MacMahon conj intro} is, however, only true modulo $t^6$ as long as $n\geqslant 4$ \cite{Atkin}. Computing the numbers $p_d^n$ for $n\geqslant 4$ is a major unsolved  problem in Combinatorics; quoting the words of Stanley, ``almost nothing significant is known'' \cite[Sec.~7.20, p.~365]{RPStanley}. While it was recently proved that the right-hand side of \eqref{eqn: MacMahon conj intro} enumerates higher dimensional partitions of fixed \emph{corner-hook volume}  \cite{Amanov-Yeliussizov}, a closed formula for its left-hand side is currently out of reach. See \Cref{sec:MacMahon-discrepancy} for more on this matter.

It is worth stressing that, despite the validity of Euler's formula, no closed (nonrecursive) formula is known for $p_d^2$, let alone for $p_d^n$. This problem originally motivated Hardy and Ramanujan, over a century ago, to look for asymptotic formulas for $p_d^2$. In \cite{Hardy-Ramanujan1} they proved that $\log\,p_d^2\approx \pi\sqrt{2d/3}$. A short while later, they improved this estimate with analytic methods, see \cite{Hardy-Ramanujan2,Hardy-Ramanujan3}. In more recent years it was confirmed  that the asymptotic growth of $p_{d}^n$ is exponentially fast for higher $n$ as well \cite{BPA_asym, Ogan_bound, PZ_proceed}; more precisely, one has
\begin{align*}
    \log\,p_d^n\approx d^{\frac{n-1}{n}}.
\end{align*}
This makes the direct computation of $p_d^n$ cumbersome and, in practice, unfeasible; for instance, the exact values of $p_d^4$ are only known for $d\leqslant 72$, and for a general $n\geqslant 1$ the numbers $p_d^n$ were computed for $d\leqslant 26$ \cite{ThePartitionsProjectWebsite, GOVINDARAJAN2013600}. 

In this paper we are interested in the \emph{exact values} of the numbers $p^\bullet_\bullet$, rather than their asymptotics. Motivated by Govindarajan's work \cite{GOVINDARAJAN2013600} and by our previous paper \cite{MOTIVEMOTIVEMOTIVE}, we study the distribution of the numbers $p_d^\bullet$ where we fix the size $d$ and allow the dimension to grow indefinitely. This means that, instead of working with the series \eqref{eqn: MacMahon conj intro}, we switch perspective and deal with the (rational) generating function
\begin{equation}
\label{eqn:H_d(t)}
\mathsf H_d(t) = \sum_{n\geqslant 0}p^{n+1}_{d}t^n.
\end{equation}
Our main results exploit a \emph{socle reduction} technique (cf.~\Cref{sec:socle-red}), through which we explicitly compute the generating functions counting higher dimensional partitions satisfying suitable constraints, \emph{refining} the original counting  $p_d^n$. We carefully analyse the rationality properties of these infinite series. The key insight of our analysis is that these refinements and their `internal' relationships, that we pinpoint throughout the paper, are enough to reconstruct all integers $p_d^\bullet$ for fixed $d$. As a biproduct, we obtain a closed expression for $\mathsf H_d(t)$ for $ d\leqslant 30$. With such expression at hand, it is immediate to compute, for instance,
\[
p_{30}^{666} = 5390806817913544023450455014935417834529246670018145780.
\]

\subsection{Refined partition counting: socle reduction technique}
\label{sec:socle-red}
As already mentioned, no closed formula is currently known for the generating function of the numbers $p_{\bullet}^n$ for $n\geqslant 4$. On the other hand, if we fix the size $d\geqslant 1$, an immediate corollary of the main result in our previous work \cite{MOTIVEMOTIVEMOTIVE} states that the generating function $\mathsf H_d(t)$ defined in \eqref{eqn:H_d(t)} is a rational function, satisfying the relation
\begin{align}\label{eqn: intro blaa}
   \mathsf H_d(t) = \frac{\mathsf h_d(t)}{(1-t)^d}\in \BZ(t),
\end{align}
where $\mathsf h_d(t)=1$ for $d=1,2,3$, and is a polynomial of degree at most $d-2$ if $d>3$. The rationality is related to the crucial observation\footnote{While this observation was already classically known \cite{Atkin}, Equation \eqref{eqn: intro blaa} follows directly from our more general result on \emph{motivic invariants} \cite{MOTIVEMOTIVEMOTIVE}.}
that the counting problem  for \emph{all} higher dimensional partitions is actually redundant, as suggested by the relation
\begin{align}\label{eqn: intro Y reduction}
      p_d^n= \sum_{k=0}^{d-1}\binom{n}{k} y_{d}^k,
\end{align}
where $y_{d}^k$ is the number of $(k-1)$-dimensional partitions of precisely \emph{embedding dimension} $k$ and size $d$ (cf.~\Cref{sec: higher dim part}). In particular, for a fixed $d\geqslant 1$, the \emph{infinitely many} numbers $p_d^\bullet$ can be computed knowing only $d-1$ initial values.

\subsubsection{Exact embedding dimension}\label{sec: sub intro 1} 
Our first main result is on the generating function
\[
\mathsf Y_e(t) = \sum_{k\geqslant 0 }y_{k+2+e}^{k+1}t^k,
\]
which should be regarded as encoding easier \emph{fundamental blocks} determining the numbers $p_d^n$.

\begin{thm}[\Cref{thm:geny}]
\label{thm: intro Y}
Fix $e \in \BZ_{>0}$. Then
\[
\mathsf Y_e(t) =\frac{ \mathsf{y}_e(t)}{(1-t)^{2e+1}}
\]
in $\BZ(t)$, where  
\[
\mathsf{y}_e(t) = \sum_{h=0}^{2e-1} \gamma_h^{(e)}t^h\in \BZ[t]
\] 
is a polynomial of  degree at most $2e-1$, whose coefficients are given by
\[
\gamma_h^{(e)} = \sum_{j=0}^h\,(-1)^{h+j}\binom{2e+1}{h-j}y_{j+e+2}^{j+1}
\]
for $0\leqslant h \leqslant 2e-1$. In particular, $\mathsf Y_e(t)$ is a rational function.
\end{thm}

We have computed the polynomials $\mathsf{y}_e(t)$ for $e\leqslant 18$ and we have listed them in \ancillary. We shall pinpoint in \Cref{cor: y determinati dai primi} how exactly the infinite sequence $(y_{k+e+1}^k)_{k\geqslant 0}$ is determined by its first $2e+1$ values.

It is worth mentioning that \cite{Atkin} introduced some  refined numbers $\overline{y}_\bullet^\bullet$, formally extracted by MacMahon's conjectural formula \eqref{eqn: MacMahon conj intro} analogously to \eqref{eqn: intro Y reduction} (cf.~\Cref{sec:conjectures}) and conjectured the positivity of the \emph{discrepancy} $e^k_d=\overline{y}_{d}^k-y_{d}^k$. We propose a conjectural rational formula 
for the generating series of the numbers $\overline{y}_{\bullet}^\bullet $ (cf.~\Cref{conj: Andrews}) and prove that implies the positivity of the discrepancy $e^k_d$ for  $d\gg 0$ (cf.~\Cref{prop:definite-positivity}).

\subsubsection{No socle in degree 1}
\label{sec: sub intro 2}
Similarly to \eqref{eqn: intro Y reduction}, it was already noticed in \cite{GOVINDARAJAN2013600, MOTIVEMOTIVEMOTIVE} that the numbers $y^\bullet_\bullet$ admit further refinements via
\begin{align*}
y_{k+e+1}^k=\sum_{j=0}^{2e} \binom{k}{j}c_e^j,
\end{align*}
where $c_e^j$ is the number of partitions of size $1+j+e$, embedding dimension $j$ and no \emph{socle} in degree 1, where the socle of a partition is the collection of its maximal elements according to the standard poset structure (see \Cref{def:socle} for the precise definition).

Define, for $x \geqslant 0$, the generating function 
\[
\mathsf C_x(t) =  
\begin{cases}
\displaystyle\sum_{e\geqslant \left\lceil \frac{x}{2}\right\rceil}\, {c_{e+1}^{2(e+1)-x}} t^{e-\left\lceil \frac{x}{2}\right\rceil} & \mbox{ if }x \mbox{ is even}\\\\
\displaystyle\sum_{e\geqslant \left\lceil \frac{x}{2}\right\rceil}\, {c_{e}^{2e-x}} t^{e-\left\lceil \frac{x}{2}\right\rceil} & \mbox{ if }x\mbox{ is odd}
\end{cases}
\]  
in $\BZ\llbracket t \rrbracket$.
Our second result is that $\mathsf C_x(t)$  becomes a nice analytic function after applying the \emph{Borel resummation} operator $\Bor \colon \BQ\llbracket t \rrbracket \to \BQ\llbracket t \rrbracket$ sending $\sum_{k\geqslant 0} f_k t^k \mapsto \sum_{k\geqslant 0} \frac{1}{k!}f_k t^k$ (cf.~\Cref{sec: Borel}).
\begin{thm}[\Cref{thm:genc}]\label{thm: intro C}
Fix $x\in\BZ_{\geqslant 0}$. Then
\begin{equation*}
   \Bor\left( \mathsf C_x(t) \right)=\frac{ \mathsf{c}_x(t)}{(1-2t)^{\frac{3}{2}+2x-\left\lceil\frac{x}{2}\right\rceil}}
\end{equation*} 
in $\BQ\llbracket t\rrbracket$, where $\mathsf c_x(t)\in \BQ[t]$ is a polynomial of  degree at most $2x-\left\lceil\frac{x}{2}\right\rceil$. In particular, $\Bor\left( \mathsf C_x(t) \right)$ is an analytic function in $t=0$.
\end{thm} 

See \Cref{thm:genc} for an explicit formula for $\mathsf c_x(t)$. We have computed these polynomials for $x\leqslant 26$ and we have listed them in \ancillary. We remark that, for $x\leqslant 5$, our result recovers the polynomials displayed in 
\cite[Sec.~4.3]{GOVINDARAJAN2013600}. 
The first polynomial that could not be extracted from the data provided in \cite{ThePartitionsProjectWebsite} is
\begin{multline*}
\mathsf c_6(t)=-\frac{4081}{6}t^9 + \frac{8491}{2}t^8 - \frac{25511}{2}t^7 + \frac{116795}{6}t^6 + \frac{1929}{4}t^5 + \frac{292709}{4}t^4 \\
+ 59667t^3 + 6773t^2 + 475t + 11.
\end{multline*}
We shall pinpoint in \Cref{cor: y determinati dai primi} how exactly the infinite sequence $(c_e^{2e-x})_{e\geqslant\lceil\frac{x}{2}\rceil}$ is determined by its first $2x-\left\lceil\frac{x}{2}\right\rceil$ values.

\subsubsection{No socle in degree 2}
Motivated by the yoga of \emph{socle reduction} outlined in \ref{sec: sub intro 2}, we show that the numbers $c_{\bullet}^\bullet $ can be further refined by new numbers $\alpha_{\bullet, \bullet}^\bullet$ and  $c_{\bullet, \bullet}^\bullet$, see Definitions \ref{def:Mpartition} and \ref{def: c molti indici}.

\begin{prop}[{\Cref{prop:recur}}]\label{prop:recurIntro}\label{prop: intro alpha}
The numbers $c_{e,(k,q,m)}^a$ satisfy the relations
\begin{equation}
    \begin{aligned}
c_{e,(k,q,m)}^0 &= 0, \\
c_{e,(k,q,m)}^a   &= \alpha_{q,m}^k\binom{a}{k} \binom{ \binom{a+1}{2}-q}{e-q-m} - \sum_{i=1}^{a-1}{ \binom{a}{i}c_{e,(k,q,m)}^i  } . 
\end{aligned}
\end{equation}
Moreover, there is an identity
\[
c_e^a = \sum_{\substack{k \leqslant a \\ q+m\leqslant  e}} c_{e,(k,q,m)}^a.
\]
\end{prop}
Thanks to the recursive structure of \Cref{prop: intro alpha}, the enumeration of higher dimensional partitions is reduced to the computation of the numbers $ \alpha_{\bullet, \bullet}^\bullet$, which count \emph{M-partitions} (cf.~\Cref{def:Mpartition}). These are higher dimensional partitions with fixed embedding dimension and number of degree 2 boxes, and admitting \emph{no socle} in degree 1 and 2.

Excluding from the counting all the higher dimensional partitions with socle in degree 2 dramatically reduces the computational complexity of $p^{\bullet}_\bullet$. In fact, M-partitions enjoy particularly nice features; for instance, their embedding dimension is always bounded by their number of quadrics, cf.~\Cref{prop:minq}. This allowed us to find exact formulas for their generating function in many cases, as we outline in \Cref{sec: exact intro}. 

In \Cref{subsec:algorithms}, we devise an algorithmic  strategy for computing the numbers $\alpha^\bullet_{\bullet, \bullet}$ in the general case, which we have implemented  in  Sage \cite{sagemath}, Mathematica \cite{Mathematica} and Macaulay2 \cite{M2} to collect a large amount of data. 

Combining the recursion of \Cref{prop: intro alpha} with the  main Theorems \ref{thm: intro Y} and \ref{thm: intro C}, the exact generating functions in \Cref{sec:exactformulas} and the algorithmic strategy presented in  \Cref{subsec:algorithms}, we were able to compute the number $p_d^n$ of higher dimensional dimensional partitions for all $n\geqslant 1$ and $d\leqslant 30$. We remark that already the number $p^n_{27}$ would be computationally out of reach for $n\geqslant 7$ using only the algorithms currently present in the literature, see e.g.~\cite[Sec.~4.3]{GOVINDARAJAN2013600}.

\smallbreak
It is worth mentioning that our socle reduction principle could ideally be extended to reduce the general partition enumeration problem to the enumeration problem of partitions with arbitrary lower bound for the socle degree. However, already in the case of constraining the socle to be in degree at least 4, finding a recursive structure analogous to the one of \Cref{prop: intro alpha} is a highly nontrivial problem, which we leave to future exploration.

\subsection{Counting M-partitions}
\label{sec: exact intro}
We introduce in \Cref{def: hydral} a special kind of M-partitions, that we call \emph{hydral partitions}. These are $(n-1)$-dimensional partitions $\lambda$ with socle in degree at least 3, and such that $h_\lambda(1)=h_\lambda(2)=n$, i.e. with minimal number of quadrics. They are counted by the integers $\alpha^n_{n,\bullet}$. The following theorem describes their count in arbitrary dimension, see also Propositions \ref{prop:hydral-n=2} and \ref{prop:dim3any} for a direct calculation of their generating function in the case $n=2$ and $n=3$ respectively, and \Cref{ex:dim4} for the case $n=4$. 

\begin{thm}[\Cref{thm: hydral structure,cor: count hydral}]
\label{thm: intro hydral}   
Fix some integer $d=1+2n+m$. Let $\alpha_{n,m}^n$ be the number   
of hydral partitions of size $d$ and embedding dimension $n$. Then
\begin{equation}
\label{eqn:hydral-count}
\sum_{m\geqslant 0} \alpha_{n,m}^nt^m=\sum_{\lambda\vdash n}\sum_{i=0}^{\mathsf t(\lambda)}\binom{n}{3i}\frac{(3i)!}{6^ii!} \mathsf f(\mathsf s^i(\lambda)) \Psi_{\mathsf{s}^i(\lambda)}(t)t^i,
\end{equation}
where $\mathsf{t,f,s}$  are auxiliary functions that we define in \Cref{not-for-thm-C} and $ \Psi_{\mathsf{s}^i(\lambda)}(t)$ is defined and proved to be rational in \Cref{lemma:techfinal}. In particular, the left-hand side of \eqref{eqn:hydral-count} is a rational function.
\end{thm} 
\Cref{thm: intro hydral} implies that it is enough to perform the algorithmic procedure explained in \Cref{subsec:algorithms} only for M-partitions with nonminimal number of quadrics. Moreover, the rationality statement in \Cref{thm: intro hydral} should be considered `optimal', in the sense that counting M-partitions with nonminimal number of quadrics, i.e.~with $h_\lambda(2)=h_\lambda(1)+1 $,  yields manifestly nonrational generating functions in any dimension, cf.~\Cref{rem:optimality}.

In a similar fashion, we prove a closed exact formula for the generating series of a special class of hydral partitions, summing over all dimensions.
\begin{thm}[\Cref{thm:expexp}]\label{thm:expexp_intro}
Fix an integer $s\geqslant 3$ and define the Hilbert--Samuel function 
\[
h_{n,s}=(1,\underbrace{n,n,\ldots ,n}_{s}).
\]
Then  $\alpha_{h_{n,s}}$ is independent of $s$ and  there is an identity
\[
\Bor\left( \sum_{n\geqslant 0} \alpha_{h_{n,s}} t^n\right)=   e^{te^t}.
\]
\end{thm}
Analogously to \Cref{thm:expexp_intro}, we prove a series of closed exact formulas counting families of \emph{(anti-)compressed} M-partitions, cf.~\Cref{prop:max1}, \ref{prop:max2}, \ref{prop:max3}. Once more, each of these exact counting eventually simplies the complexity of the original counting problem $p_\bullet^\bullet$.
\subsection{Geometry and Physics} 
Fix integers $n>0$ and $d \geqslant 0$. Let $Z$ be a smooth $n$-dimensional complex variety, $E$ a locally free sheaf of rank $r\geqslant 1$ over $Z$. Form the \emph{Quot scheme of points} $\Quot_Z(E,d)$, namely the fine moduli space of coherent sheaf quotients $E \onto Q$, where $\chi(Q) = d$ \cite{Grothendieck_Quot}. When $r=1$, this space reduces to the \emph{Hilbert scheme of points} $\Hilb^d(Z)$, which is almost always singular (precisely when $n\geqslant 3$ and $d\geqslant 4$) and pathological in a precise sense \cite{Jelisiejew-pathologies}. The integers $p_\bullet^n$ determine the topological Euler characteristic of $\Quot_Z(E,\bullet)$ via
\begin{equation}
\label{eqn:quot-gf}
\sum_{d \geqslant 0} \chi(\Quot_Z(E,d))t^d = \left(\sum_{d\geqslant 0}p_d^nt^d\right)^{r\chi(Z)},
\end{equation}
independently on $E$, which should be seen as a numerical shadow of the \emph{motivic} identities proved in \cite{ricolfi2019motive,GLMHilb}. We remark that, through the viewpoint of algebraic geometry, the failure of the motivic analagoue of  MacMahon's conjecture is deeply related with the singularities of the Hilbert scheme  $\Hilb^d(\BA^n)$, cf.~\cite{MOTIVEMOTIVEMOTIVE}.

The connection between the enumeration of higher dimensional partitions and Quot schemes is  ubiquitous in enumerative geometry (Donaldson--Thomas theory),  geometric representation theory and string theory, where the main focus is on the cases $n=3,4$, see \cite{MNOP1,BBS,DavisonR,BR18,cazzaniga2020higher,Cazzaniga:2020aa,doi:10.1142/S0219199724500196,Virtual_Quot,FMR_K-DT,MR_nested_Quot,double-nested-1,CKM_K_theoretic, Mon_canonical_vertex, fasola2024tetrahedroninstantonsdonaldsonthomastheory, nekrasov2017magnificent, Nekrasov_M-theory,ORV_classical_crystals} for a nonexhaustive list of references.

\subsection*{Acknowledgements}
We would like to thank Geoffrey Critzer, Alessio D'Alì, Paolo Lella and Alessio Sammartano for useful discussions. Special thanks to Giorgio Gubbiotti for generously sharing his expertise and technical support on some rationality questions. Finally, we thank Andrea Sangiovanni for suggesting the link with hypergeometric functions (\Cref{rmk:hypergeom}). 

S.M. is supported by the FNS Project 200021-196960 ``Arithmetic aspects of moduli spaces on
curves''.  R.M. is partially supported by the PRIN project 2022L34E7W ``Moduli Spaces and Birational Geometry''. A.R. is partially supported by the PRIN project 2022BTA242 ``Geometry of algebraic structures: moduli, invariants, deformations''. All authors are members of the GNSAGA - INdAM.

\section{Background and notation}
\subsection{Conventions}
We set $\BN = \BZ_{\geqslant 0}$. For $n\in \BZ$ and $k\in \BN$, we set
\[
\binom{n}{k}=  \frac{\prod_{i=0}^{k-1}(n-i)}{k!},
\]
which coincides with the classical binomial coefficient if $n\geqslant 0$.
Moreover, for $n\in \BZ_{\geqslant 1}$, the lattice $\BN^n$ is endowed with its stanfard component-wise poset structure throughout. 

The cardinality of a finite set $S$ will be denoted $\lvert S \rvert$.
    
\subsection{Higher dimensional partitions}\label{sec: higher dim part}
We introduce the main object of study of this paper.
\begin{definition}
\label{def:partition}
Fix $n,d \in \BZ_{\geqslant 0}$. An \emph{$(n-1)$-dimensional partition} is a collection of $d$ points $\lambda =\set{\ba_1,\ldots,\ba_d}\subset \BN^{n}$ such that if  $\by \in \BN^{n}$ satisfies $\by\le\ba_i$ for some $i=1,\ldots,d$, then $\by \in \lambda$. We call $|\lambda|=d$ the \emph{size} of $\lambda$ and we denote by $\mathrm{P}^n_d$ the set of $(n-1)$-dimensional partitions of size $d$, and by $p_d^n$ the cardinality $\lvert \mathrm{P}^n_d\rvert$.
\end{definition}

A partition $\lambda \in \mathrm{P}^n_d$ can alternatively be visualised as a stack of boxes in $\BR^{n-1}$, labelled by nonincreasing positive integers, whence the name \emph{$(n-1)$-dimensional partitions}.

For instance,
\[
p_{d}^n = 
\begin{cases}
    1 & \mbox{if } d\in\set{0,1},\,n=0 \\
    0 & \mbox{if }d>1,\,n=0 \\
    1 & \mbox{if } n=1 \\
    1 & \mbox{if }d=0. 
\end{cases}
\]
Elements of $\mathrm{P}^2_d$ (resp.~$\mathrm{P}^3_d$) are known in the literature as \emph{linear} (resp.~\emph{plane}) partitions. Linear partitions of size $d>0$ are in natural bijections with Young diagrams of size $d$. 

\smallbreak
Next, we recall a few constructions and invariants attached to higher dimensional partitions.

Given a lattice point $\by=(y_1,\ldots,y_n)\in \BN^n$, we define its \emph{degree} by
\[
\deg(\by)=\sum_{i=1}^ny_i.
\]
Let $\lambda \in \mathrm{P}^n_d$ be a partition. For $k\geqslant 1$, we set
\[
\lambda_{\geqslant k} =\Set{\by\in \lambda |\deg(\by)\geqslant k }\subset\lambda,
\]
and similarly we denote by $\lambda_{=k}$ (resp.~$\lambda_{\leqslant k}$) the set of lattice points of degree exactly $k$ (resp.~at most $k$). The \emph{Hilbert--Samuel function} of $\lambda$ is the function 
\[
\begin{tikzcd}[row sep=tiny]
 \BN\arrow[r,"h_\lambda"]&\BN\\
    i\arrow[r,mapsto] & {\lvert }\lambda_{=i}{\rvert},
\end{tikzcd}
\]
which we sometimes represent as a tuple of positive integers $h_\lambda=(h_0,h_1, \dots, h_{\ell(\lambda)})$, calling $\ell(\lambda)$ the \emph{length} of $\lambda$. We also set $\lvert h_\lambda \rvert = \lvert \lambda \rvert$. Note that the length of $\lambda$ satisfies
\[
\ell(\lambda)=\max\Set{\deg(\by)|\by\in\lambda  }.
\]
We call $h_\lambda(1)$ the \emph{embedding dimension} of $\lambda$. 
The \emph{socle} $\Soc(\lambda)$ is defined to be the set of maximal elements in $\lambda$, namely
\begin{equation}
\label{def:socle}
\Soc(\lambda)=\Set{\by\in\lambda| \by \not< \ba \mbox{ for all } \ba \in \lambda}\subset \lambda.
\end{equation}
A basic, but key observation is that by the very definition of partition, the socle $\Soc(\lambda)$ determines $\lambda$ uniquely.

Finally, the \emph{socle type} of $\lambda$ is the function  
\[
\begin{tikzcd}[row sep=tiny]
\BN\arrow[r,"e_\lambda"]&\BN\\
    i\arrow[r,mapsto] & \lvert \Soc(\lambda) \cap \lambda_{=i}\rvert,
\end{tikzcd}
\]
which we sometimes represent as a tuple of integers $e_\lambda=(e_0,e_1, \dots, e_{\ell(\lambda)})$.

\begin{notation}
\label{notation:linear-partitions}
Given a linear partition $\lambda \in \mathrm{P}^2_d$, there are two different (equivalent) ways to represent $\lambda$, namely
\begin{itemize}
    \item [\mylabel{linear-partition-1}{(1)}] $\lambda = (\lambda_1,\lambda_2,\ldots,\lambda_\ell)$ where $\lambda_1\geqslant \lambda_2\geqslant \cdots\geqslant \lambda_\ell>0$ and $d = \sum_j\lambda_j$, and
    \item [\mylabel{linear-partition-2}{(2)}] $\lambda = (1^{\alpha_1}2^{\alpha_2}\cdots d^{\alpha_d})$ where $d = \sum_ii\alpha_i$ (and some $\alpha_i$ might be 0, in which case they are omitted).
\end{itemize}
In both cases, we write `$\lambda \vdash d$' to say that $\lambda \in \mathrm P_d^2$. The notation in \ref{linear-partition-2} means that $\lambda$ consists of $\alpha_i$ parts of size $i$ for all $i=1,\ldots,d$, so that $\ell = \sum_i\alpha_i$ counts the number of columns in the Young diagram attached to $\lambda$. To pass from \ref{linear-partition-1} to \ref{linear-partition-2}, consider the string $\lambda' = (\lambda_1',\lambda_2',\ldots,\lambda_{s}')$ of decreasing \emph{distinct} parts of $\lambda$ (so that $\lambda_1=\lambda_1'>\lambda_2'>\cdots>\lambda_s'$ and $s\leqslant \ell$). Each $\lambda_j'$ satisfies $1\leqslant \lambda_j'\leqslant d$ and carries a positive multiplicity $\alpha_{\lambda_j'}\leqslant d$, counting the number of times it occurred in $(\lambda_1,\lambda_2,\ldots,\lambda_\ell)$. For example, the Young diagram
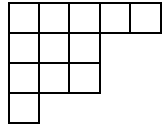
\[
\begin{matrix}
       \yng(5,3,3,1) \quad 
       \ytableausetup{boxsize=1.09em}
       \ytableausetup{boxframe=0.02em}\ytableausetup{aligntableaux=bottom}
   \end{matrix}
\]
corresponds to $\lambda = (4,3,3,1,1) = (1^2 3^2 4^1)$, has size $12$, satisfies $\lambda'=(4,3,1)$ and has $2+2+1 = 5$ parts.

The \emph{automorphism group} of $\lambda \in \mathrm{P}^2_d$ is, by definition, the group $\aut(\lambda) = \prod_i\FS_{\alpha_i}$, where $\FS_a$ is the symmetric group on $a$ letters. In the above example, the automorphism group is $\FS_2\times\FS_2$.
\end{notation}

\subsection{Apolarity}
Let $S\subset \BN^n$ be any finite subset. The \emph{apolar partition} $S^\perp$ associated to $S$ is the partition  
\[
S^\perp = \Set{\by\in\BN^n|\by\le\ba \mbox{ for some } \ba\in S} \subset \BN^n. 
\]
We say that a finite subset $S\subset \BN^n$ is \emph{admissible} if any two elements in $S$ are not comparable. By construction, for all partitions $\lambda \in \mathrm{P}_d^n$ and finite admissible subsets $S\subset \BN^n$ we have that
\begin{align*}
     \Soc(\lambda)^\perp&=\lambda,\\
     \Soc(S^\perp)&= S. 
\end{align*}
In particular, there is a bijection between partitions and admissible subsets $S\subset \BN^n$.

\subsection{Monomial ideals}
\label{sec: monomial ideals} 
Fix $n \in \BZ_{\geqslant 0}$. There is a bijective correspondence between $(n-1)$-dimensional partitions and monomial ideals of finite colength in $\BC[x_1, \ldots, x_n]$. To a partition $\lambda\in \mathrm{P}_d^n$, one associates the monomial ideal
\[
I_\lambda=\Braket{x_1^{i_1}\cdots x_n^{i_n}|(i_1,\ldots,i_n) \in \BN^n\setminus \lambda}\subset \BC[x_1, \ldots, x_n],
\]
with quotient
\[
\BC[x_1, \dots, x_n]/I_\lambda=\bigoplus_{(i_1,\ldots,i_n) \in\lambda} \BC\cdot x_1^{i_1}\cdots x_n^{i_n}.
\]
We now reformulate the apolar construction for ideals and their quotients (see e.g.~\cite[Sec.~2.5]{MOTIVEMOTIVEMOTIVE}). For a finite subset $T\subset \BC[x_1, \dots, x_n]$, define its \emph{apolar ideal}
\[
T^\perp=\Braket{x_1^{\alpha_1}\cdots x_n^{\alpha_n}\,|\,\displaystyle\frac{\partial^{\sum_{i=1}^n \alpha_i} }{\partial^{\alpha_1} x_1\cdots \partial^{\alpha_n} x_n} p=0,  \mbox{ for all } p\in T,\ \alpha_i\in\BN     }\subset  \BC[x_1, \dots, x_n].
\]
To an admissible subset $S\subset \BN^n$, we associate the collection of monomials
\[
\widetilde{S}=\Set{x_1^{i_1}\cdots x_n^{i_n}|(i_1,\ldots,i_n) \in S}\subset \BC[x_1, \dots, x_n],
\]
whose apolar ideal $\widetilde{S}^\perp$ corresponds to the apolar partition $S^\perp$. This formulation will be a useful language to describe partitions in terms of the monomial ideals generated via apolarity by their socles.

Whenever it will be clear from the context, we will identify partitions with their associated monomial ideals.

\subsection{Borel transform}
\label{sec: Borel}
Several formulas that we shall prove in this paper are most conveniently formulated in terms of the \emph{Borel transform}, namely the operator
\[
\begin{tikzcd}[row sep =tiny]
\BQ\llbracket t \rrbracket\arrow[r,"\Bor"]&
\BQ\llbracket t \rrbracket\\
\displaystyle\sum_{k\geqslant 0} f_k t^k\arrow[r,mapsto]& \displaystyle\sum_{k\geqslant 0}\frac{f_k}{k!} t^k.
\end{tikzcd}
\]
We denote by $\IB$ its inverse. We now prove two identities that will be used throughout the paper.

\begin{lemma}
\label{lemma:borelmainprop}
For every $x \in \BZ_{>0}$ there is an identity
   \[
    \IB \left(\frac{1}{(1-2t)^{\frac{1+2x}{2}}}\right)=\sum_{n\geqslant 0 }\frac{(2n+2x-1)!!}{(2x-1)!!} t^n
   \] 
\end{lemma}

\begin{proof}
By the binomial theorem with rational exponents, we have, for every rational number $m \in \BQ$, the identity
\begin{equation}\label{Rational-Binomial}
(1+u)^m = \sum_{i\geqslant 0}\left(\frac{1}{i!}\prod_{k=0}^{i-1}\,(m-k) \right) u^i.
\end{equation}
Specialising to $u = -2t$ and $m = -x-1/2$, we obtain
\[
\frac{1}{(1-2t)^{\frac{1+2x}{2}}} = \sum_{i\geqslant 0} \left( \frac{(-2)^i}{i!}\prod_{k=0}^{i-1}\,(-x-1/2-k)\right)t^i.
\]
The claim then reduces to proving the identity
\[
(-2)^i\prod_{k=0}^{i-1}\,(-x-1/2-k) = \frac{(2i+2x-1)!!}{(2x-1)!!}
\]
for every $i \geqslant 0$. But this is clear, since
\[
(2x-1)!! (-2)^i\prod_{k=0}^{i-1}\,(-x-1/2-k) 
=  (2x-1)!!\prod_{k=0}^{i-1}\, (2x+1+2k)
=(2i+2x-1)!!
\]
as required.
\end{proof}

\begin{lemma}
\label{lemma:techborel}
For every $s\in \BZ_{>0}$ there is an identity
\[
\IB \left(\frac{ (2s-1)!!(s+t)}{(1-2t)^{\frac{3+2s}{2}}}\right)=\sum_{i\geqslant 0} (s+i)(2s+2i-1)!!t^{i}.
\]
\end{lemma}

\begin{proof}
We split the argument of $\IB$ into two functions, corresponding to the sum $s+t$. By \Cref{Rational-Binomial}, we have
\begin{align*}
\frac{s(2s-1)!!}{(1-2t)^{\frac{3+2s}{2}}} 
&= s(2s-1)!!\sum_{i\geqslant 0}\left(\frac{1}{i!}\prod_{k=0}^{i-1}(3+2s+2k) \right) t^i \\
\frac{t(2s-1)!!}{(1-2t)^{\frac{3+2s}{2}}} 
&=(2s-1)!!\sum_{i\geqslant 1}\left(\frac{1}{(i-1)!}\prod_{k=0}^{i-2}(3+2s+2k) \right) t^i.
\end{align*}
Applying the operator $\IB$, we thus obtain
\begin{equation}\label{equations-IB_s}
\begin{split}
\IB\left(\frac{s(2s-1)!!}{(1-2t)^{\frac{3+2s}{2}}} \right) 
&=s(2s-1)!!\sum_{i\geqslant 0} \left(\prod_{k=0}^{i-1}(3+2s+2k) \right) t^i \\
\IB\left(\frac{t(2s-1)!!}{(1-2t)^{\frac{3+2s}{2}}} \right)
&=(2s-1)!!\sum_{i\geqslant 1}\left(i\prod_{k=0}^{i-2}(3+2s+2k) \right) t^i.
\end{split}
\end{equation}
Let us check the coefficients of $t^0$ and $t^1$ separately before dealing with the general case. As for $t^0$, only the first function contributes, and gives $s(2s-1)!!$, as it should. As for $t^1$, we obtain
\[
s(2s-1)!!(3+2s)+(2s-1)!! = (2s-1)!!(1+3s+2s^2) = (2s-1)!!(2s+1)(s+1),
\]
which coincides with the required value $(s+1)(2s+1)!!$. For the cases $i\geqslant 2$, summing the contributions from Equations \eqref{equations-IB_s} yields
\begin{align*}
s(2s-1)!!&\prod_{k=0}^{i-1}(3+2s+2k)+i(2s-1)!!\prod_{k=0}^{i-2}(3+2s+2k) \\
&=\left((2s-1)!!\prod_{k=0}^{i-2}(3+2s+2k)\right)\left( s(2s+2i+1)+i\right) \\
&=(2s-1)!!(2s+3)(2s+5)\cdots (2s+2i-1)(s+i)(2s+1) \\
&=(2s+2i-1)!!(s+i),
\end{align*}
as required.
\end{proof}

\section{Generating functions of partitions}
In the next subsection, we introduce two distinguished collections $\mathrm{Y}^\bullet_\bullet$ and $\mathrm{C}^\bullet_\bullet$ of higher dimensional partitions, whose associated counting theories reconstruct the original counting problem for $p_{\bullet}^\bullet$, as we establish by means of the `inversion formulas' in \Cref{prop: inversion} below. In subsections \ref{sec:proof-Y-intro} and \ref{sec: proof C main} we prove Theorems \ref{thm: intro Y} and \ref{thm: intro C} respectively.

\subsection{Refined partition counting and inversion formulas}\label{subsec:ref count}

Fix $k,d,e \in \BN$. Define the sets
\begin{equation}
\label{Y-and-C}
\begin{split}
\mathrm{Y}^k_d 
&=\Set{\lambda \in \mathrm{P}^k_d| h_\lambda(1)=k}\subset \mathrm{P}^k_d, \\
\mathrm{C}^k_e 
&=\Set{\lambda \in \mathrm{Y}^k_{1+k+e} | \Soc(\lambda)\subset \lambda_{\geqslant 2}}\subset \mathrm{Y}^k_{1+k+e}, \\
\mathrm{D}^k_e &=\Set{\lambda \in \mathrm{C}^k_e | \lambda_{=2}\mbox{ is M-stable, cf.~\Cref{subsec:algorithms}}}\subset \mathrm{C}^k_e.
\end{split}
\end{equation}
The sets $\mathrm{C}_\bullet^\bullet$ are key to our \emph{socle reduction} technique. The sets $\mathrm{D}^k_e$ are only defined for $e\geqslant 0$ and $0 \leqslant k\leqslant  2e-\lceil e/2\rceil$. We will not expand on their definition as they will play only an auxiliary role in the proof of \Cref{thm: intro C}, but we refer to \cite{GOVINDARAJAN2013600} for more details on $\mathrm D_\bullet^\bullet$. We denote the cardinalities of the sets defined in \eqref{Y-and-C} by
\[
y^k_d=\lvert \mathrm{Y}^k_d\rvert, \qquad  c^k_e=\lvert \mathrm{C}^k_e\rvert, \qquad  d^k_e=\lvert \mathrm{D}^k_e\rvert.
\]

We present now a list of exact values of the numbers $y^k_d, c^k_e$ in some elementary limit cases (see also \Cref{ex:somenumbers}).

\begin{prop}
\label{rem:limitcases}
We have
\begin{itemize}
    \item [\mylabel{it:1}{\normalfont{(1)}}]  $y^k_d=0$ for $k\geqslant d$, or  $d>1$ and $k=0$,
    \item [\mylabel{it:2}{\normalfont{(2)}}]  $c^k_e=0$ for $k>2e$, or $k=0$ and $e>0$,
    \item [\mylabel{it:3}{\normalfont{(3)}}]  $c^k_e=1$ for $k=e=0$ or $k=1$ and $e>0$ or $k=2$ and $e=1$,
    \item [\mylabel{it:4}{\normalfont{(4)}}]  $c^{2e}_e= (2e-1)!!$ for $e>0$,
    \item [\mylabel{it:5}{\normalfont{(5)}}]  $c^{2e-1}_e= e (2e-1)!!$ for $e>0$.
    \item [\mylabel{it:6}{\normalfont{(6)}}]  $y_{n+1}^2=\Coef_{q^{n+1}}\left[\!\begin{smallmatrix} 2n \\ n \end{smallmatrix}\!\right]_q$,
where
\[
\left[\begin{matrix}
    a\\b
\end{matrix}\right]_q=\frac{\prod_{j=1}^a(q^j-1)}{\prod_{j=1}^{b}(q^j-1)\prod_{j=1}^{a-b}(q^j-1)}\in\BZ[q]. 
\] 
    \item [\mylabel{it:7}{\normalfont{(7)}}]  $y_{n}^2=p_{n}^2-2$.
    \item [\mylabel{it:8}{\normalfont{(8)}}] $c_{n-2}^2=p_{n+1}^2-4=y_{n+1}^2-2$.
\end{itemize}
\end{prop}

\begin{proof} Points \ref{it:1}, \ref{it:2} and \ref{it:3} follow from the definition of $y^k_d, c_e^k$. Points \ref{it:4} and \ref{it:5} follow from \cite[Section 2.5]{GOVINDARAJAN2013600}. Point \ref{it:6} follows from \ref{it:7} after noting that
\[
\Coef_{q^{n+1}}
\begin{bmatrix}
2n\\n
\end{bmatrix}_q = p^2_{n+1}-2.
\]
Indeed, it is well-known that the coefficient of $q^{n+1}$ in this $q$-binomial coefficient counts the number of partitions of $n+1$ fitting in an $n\times n$ square. Of all the partitions of $n+1$, only the curvilinear ones do not have this property.

Relation \ref{it:7} follows from the observation that the only two linear partitions having nonmaximal embedding dimension are the curvilinear ones. Similarly, for \ref{it:8}, the numbers $c_{n-2}^2$ exludes the curvilinear partitions and the partitions $(2,1,\ldots,1)$ and $(n,1)$, as these are the only ones having socle in degree 1.
\end{proof}

\begin{prop}[Inversion formulas]
\label{prop: inversion}
The following statements hold.
\begin{itemize}
\item [\mylabel{inversion-1}{\normalfont{(1)}}]
The collections of integers $p_\bullet^\bullet$ and $y_\bullet^\bullet$ determine each other via the relations
\[
p_d^n = \sum_{k=0}^{d-1}\binom{n}{k} y_{d}^k,\qquad 
y^n_d = \sum_{j=0}^{n} (-1)^{n+j}\binom{n}{j} p^j_d,
\]
for each $n\geqslant 0$ and $d\geqslant 1$.
\item [\mylabel{inversion-2}{\normalfont{(2)}}]
The collections of integers $y_\bullet^\bullet$ and $c_\bullet^\bullet$ determine each other via the relations
\[
y_{k+e+1}^k = \sum_{x=0}^{2e} \binom{k}{x}c_e^x,\qquad
c^k_e = \sum_{j=0}^{k} (-1)^{k+j}\binom{k}{j} y^j_{e+j+1},
\]
for each $k,e \geqslant 0$.
\item [\mylabel{inversion-3}{\normalfont{(3)}}] 
The collections of integers $c_\bullet^\bullet$ and $d_\bullet^\bullet$ determine each other via the relations
\begin{equation}
\label{inversion-cd}
\begin{split}
c_e^{2e-x}&=\sum^{\min(e,2x)}_{y=\left\lceil \frac{x}{2}\right\rceil} \frac{(2e-x)!}{( 2e-2y)!! (2y -x )!} d_{y}^{2y-x}, \\
d^{x}_e &= \sum ^{\left\lfloor \frac{x}{2} \right\rfloor}_{y=0}(-1)^{y}\frac{x!}{(2y)!!(x-2y)! }c_{e-y}^{x-2y},
\end{split}
\end{equation}
for each $e\geqslant 0$ and $0 \leqslant x\leqslant  2e-\lceil e/2\rceil$.
\end{itemize}
\end{prop}

\begin{proof}
The first identity in \ref{inversion-1} follows from the definitions,\footnote{Alternatively, use \cite[Thm.~3.4, Eqn.~(3.1)]{MOTIVEMOTIVEMOTIVE}, and set $\BL=1$.} and so does the first identity in \ref{inversion-2}. The second identity in \ref{inversion-1} is a direct check using the first identity.\footnote{Alternatively, use \cite[Lemma 5.1]{MOTIVEMOTIVEMOTIVE}, and set $\BL=1$.} The same happens for the second identity in \ref{inversion-2}.

As for \ref{inversion-3}, we know from \cite[Sec.~2.6]{GOVINDARAJAN2013600} that the first relation in \eqref{inversion-cd} is satisfied. It then remains to confirm that the two equations in \eqref{inversion-cd} compose to the identity. This is indeed the case, as for $\left\lceil\frac{x}{2}\right\rceil\leqslant e < 2x$, we have
\begin{align*}
    c_{e}^{2e-x} & = \sum^{e}_{  y=\left\lceil \frac{x}{2}\right\rceil} \frac{(2e-x)!}{( 2e-2y)!! (2y -x )!} d_{y}^{2y-x}\\
    & = \sum^{e}_{  y=\left\lceil \frac{x}{2}\right\rceil} \sum ^{y-\left\lceil \frac{x}{2} \right\rceil}_{z=0} \frac{(-1)^{z}(2e-x)! }{( 2e-2y)!!(2z)!!(2y-x-2z)! }  c_{y-z}^{2y-x-2z}\\
    & = \sum^{e}_{  y=\left\lceil \frac{x}{2}\right\rceil} \sum ^{y}_{z=\left\lceil \frac{x}{2} \right\rceil} \frac{(-1)^{z+y}(2e-x)! }{( 2e-2y)!!(2y-2z)!!(2z-x)! }  c_{z}^{2z-x}\\
    & = \sum^{e}_{  y=\left\lceil \frac{x}{2}\right\rceil} \sum ^{y}_{z=\left\lceil \frac{x}{2} \right\rceil}(-1)^{z+y}\binom{2e-x}{2z-x} \frac{(2e-2z)! }{( 2e-2y)!!(2y-2z)!!  }  c_{z}^{2z-x}\\
    & = \sum^{e}_{  y=\left\lceil \frac{x}{2}\right\rceil} \sum ^{y}_{z=\left\lceil \frac{x}{2} \right\rceil}(-1)^{z+y}\binom{2e-x}{2z-x}\binom{e-z}{e-y} \frac{(2e-2z)!}{(2e-2z)!!}  c_{z}^{2z-x}\\
    & = \sum^{e}_{  z=\left\lceil \frac{x}{2}\right\rceil} (-1)^{z}\binom{2e-x}{2z-x}\frac{(2e-2z)!}{(2e-2z)!!}  c_{z}^{2z-x}\sum ^{e}_{y=z}(-1)^{y}\binom{e-z}{e-y} \\
    & = \sum_{z=\left\lceil\frac{x}{2}\right\rceil}^{e}(-1)^z\binom{2e-x}{2z-x}\binom{ e-z-1}{2x-z}\frac{(2e-2z)!}{(2e-2z)!!}c_{z}^{2z-x}.
\end{align*}
To conclude, notice that 
\[
(-1)^z\binom{2e-x}{2z-x}\binom{ e-z-1}{2x-z}\frac{(2e-2z)!}{(2e-2z)!!}=\delta_{e,z},
\]
where $\delta_{e,z}$ is the Kronecker symbol. 
\end{proof}

As a corollary, we obtain a recursive expression for the sequences  $(y_{k+e+1}^k)_{k\geqslant 0}$ and $ (c_{e}^{2e-x})_{e\geqslant 0}$ in terms of finitely many initial values.

\begin{corollary}
\label{cor: y determinati dai primi}
 
The sequence $(y_{k+e+1}^k)_{k\geqslant 0}$ is determined by its first $2e+1 $ values. More precisely, for $k>2e$ we have
\begin{align*}
y_{k+e+1}^k&=\sum_{j= 0}^{2e} (-1)^j\binom{k}{j}\binom{k-j-1}{2e-j}y_{j+e+1}^j.
\end{align*} 
Similarly, the sequence $(c_{e}^{2e-x})_{e\geqslant 0}$ is determined by its first $2x-\left\lceil\frac{x}{2}\right\rceil $ values. More precisely, for $e> 2x$ we have
\[
c_{e}^{2e-x}= \sum_{z=\left\lceil\frac{x}{2}\right\rceil}^{2x}(-1)^z\binom{2e-x}{2z-x}\binom{ e-z-1}{2x-z}(2e-2z-1)!!c_{z}^{2z-x}
\]
\end{corollary}

\begin{proof}
The first equation follows from \Cref{prop: inversion} and a direct calculation. The proof of the second one is identical to that of  \Cref{prop: inversion}.
\end{proof}

\Cref{cor: y determinati dai primi} recovers (and extends) the formulas in \cite[Sec.~4.3]{GOVINDARAJAN2013600}. 

\subsection{Proof of \texorpdfstring{\Cref{thm: intro Y}}{}} 
\label{sec:proof-Y-intro}
In this section we prove \Cref{thm: intro Y} from the Introduction. Define the generating function
\[
\mathsf Y_e(t) = \sum_{k\geqslant 0}\,y_{k+2+e}^{k+1}t^k \,\in\,\BZ\llbracket t \rrbracket.
\] 

\begin{theorem}\label{thm:geny}
Fix $e \in \BZ_{>0}$. Then
\[
\mathsf Y_e(t) =\frac{ \mathsf{y}_e(t)}{(1-t)^{2e+1}}
\]
in $\BZ(t)$, where  
\[
\mathsf{y}_e(t) = \sum_{h=0}^{2e-1} \gamma_h^{(e)}t^h\in \BZ[t]
\] 
is a polynomial of  degree at most $2e-1$, whose coefficients are given by
\[
\gamma_h^{(e)} = \sum_{j=0}^h\,(-1)^{h+j}\binom{2e+1}{h-j}y_{j+e+2}^{j+1}
\]
for $0\leqslant h \leqslant 2e-1$. In particular, $\mathsf Y_e(t)$ is a rational function.
\end{theorem} 

\begin{proof}
   By definition we have
    \begin{align*} 
\mathsf Y_e(t) 
&= \sum_{k\geqslant 0}\,y_{k+2+e}^{k+1}t^k \\
&= \sum_{k\geqslant 0}\,\left(\sum_{x=0}^{2e} \binom{k+1}{x}c_{e}^x \right)t^k & \mbox{ by \Cref{prop: inversion}}  \\   
&
=\sum_{x=0}^{2e} c_{e}^x\frac{t^{x-1}}{(1-t)^{x+1}}&   \\
&=\frac{1}{(1-t)^{2e+1}}\sum_{x=0}^{2e} c_{e}^x(1-t)^{2e-x}t^{x-1}&   \\   
&=\frac{1}{(1-t)^{2e+1}}\sum_{x=0}^{2e}\sum_{h=x}^{2e}\binom{2e-x}{h-x}(-1)^{h+x} c_{e}^xt^{h-1}&   \\   
&=\frac{1}{(1-t)^{2e+1}}\sum_{h=1}^{2e}\sum_{x=1}^{h}\binom{2e-x}{h-x}(-1)^{h+x} c_{e}^xt^{h-1}. 
\end{align*}  
This proves the rationality claim. We now check that 
\[
\sum_{h=1}^{2e}\sum_{x=1}^{h}\binom{2e-x}{h-x}(-1)^{h+x} c_{e}^xt^{h-1} = \mathsf y_e(t),
\]
with $\mathsf y_e(t)$ as in the statement. The coefficient of $t^h$ in the left-hand side is  
\begin{align*}
\sum_{x=0}^{h}(-1)^{h+x}  \binom{2e-x-1}{h-x}c_{e}^{x+1}
&=  \sum_{x=0}^{h}(-1)^{h+x}  \binom{2e-x-1}{h-x}\left( \sum_{j=1}^{x+1} (-1)^{x+1+j}\binom{x+1}{j} y^j_{e+j+1}\right)\\      
&=  \sum_{x=0}^{h} \sum_{j=0}^{x} (-1)^{h+j}\binom{2e-x-1}{h-x}   \binom{x+1}{j+1} y^{j+1}_{e+j+2} \\   
&=   \sum_{j=0}^{h} (-1)^{h+j}\binom{2e+1}{h-j} y^{j+1}_{e+j+2},     
\end{align*}
where the first identity follows from \Cref{prop: inversion}\ref{inversion-2} and \Cref{rem:limitcases}\ref{it:1}, and the last identity follows from \cite[Cor.~4.5]{MOTIVEMOTIVEMOTIVE} after a suitable change of variables.
\end{proof}

We prove two relations among the coefficients of the polynomial $\mathsf y_e(t)$, which reduce the computational data needed to determine the series $\mathsf Y_e(t)$.
\begin{corollary}\label{prop:relationcoeffy}
Fix $e>0$. There are identities
\begin{itemize}
\item [\mylabel{gamma1}{\normalfont{(1)}}] $\displaystyle\sum _{0\leqslant i \leqslant 2e-1} \gamma_{i}^{(e)}=(2e-1)!!$, \\
\item [\mylabel{gamma2}{\normalfont{(2)}}] $\displaystyle\sum _{0\leqslant i \leqslant 2e-1} (i+1)   \gamma_{i}^{(e)} = e(2e-1)!!$. 
\end{itemize} 
\end{corollary} 

\begin{proof}
The proof of \ref{gamma1} is given by
    \begin{align*}
        (2e-1)!!&=c_e^{2e} & \mbox{\Cref{rem:limitcases}\ref{it:4}}\\
        &=\sum_{j=0}^{2e} (-1)^{j}\binom{2e}{j} y^j_{e+j+1} &\mbox{\Cref{prop: inversion}\ref{inversion-2}} 
        \\
        &=\sum_{j=1}^{2e}  \sum_{\alpha=0}^j(-1)^{\alpha } \binom{2e+1}{\alpha}  y^j_{e+j+1} \\ 
         &=\sum_{\beta=0}^{2e-1}  y^{\beta+1}_{e+\beta+2} \sum_{\alpha=0}^{\beta+1}   (-1)^{\alpha } \binom{2e+1}{ \alpha}  \\  
         &=\sum_{\beta=0}^{2e-1}  y^{\beta+1}_{e+\beta+2} \sum_{\alpha=\beta +2}^{2e+1}   (-1)^{\alpha +1} \binom{2e+1}{ \alpha}  \\
         &=\sum_{\beta=0}^{2e-1}  y^{\beta+1}_{e+\beta+2} \sum_{\overline{\alpha}=\beta +2}^{2e+1}   (-1)^{\overline{\alpha} +\beta} \binom{2e+1}{\overline{\alpha}-\beta-2}\\  
        & =\sum_{\beta=0}^{2e-1}   \sum_{\overline{\alpha}=\beta }^{2e-1}   (-1)^{\overline{\alpha} +\beta} \binom{2e+1}{\overline{\alpha}-\beta} y^{\beta+1}_{e+\beta+2}  \\  
         &=\sum_{\overline{\alpha}=0 }^{2e-1}\sum_{\beta=0}^{\overline{\alpha}}      (-1)^{\overline{\alpha} +\beta} \binom{2e+1}{\overline{\alpha}-\beta} y^{\beta+1}_{e+\beta+2} \\
         &= \mathsf y_e(1), &
    \end{align*}  
where in the third equality we used \cite[Lemma 4.11]{MOTIVEMOTIVEMOTIVE} with $\BL=1$ and where the final identity follows from \Cref{thm:geny}. 

The proof of \ref{gamma2} follows from an analogous formal manipulation involving binomial coefficients and \Cref{rem:limitcases}.
\end{proof}
We shall now list, for all $d\geqslant 1$, some values of $y^k_d$, which along with the relations of \Cref{prop:relationcoeffy} decrease the computational complexity of the polynomial $\mathsf{y}_e(t)$ of \Cref{thm:geny}.

\begin{example}
\label{ex:somenumbers}
Let $s_k=\binom{k+1}{2} $ be the number of linearly independent quadrics in $k$ variables. Then,  for all $d\geqslant 1$, there are identities
\begin{itemize}
\item [\mylabel{Y1}{\normalfont{(1)}}] $y_{d}^1=1, $
\item [\mylabel{Y2}{\normalfont{(2)}}] $y_{d}^2=p^2_d-2,$
\item [\mylabel{Y3}{\normalfont{(3)}}] $y_{d}^3=p^3_d-3(y_{d}^2+1),$
\item [\mylabel{Y4}{\normalfont{(4)}}] $y_{d}^{d-5}=  \binom{s_{d-5}}{4} + (d-5)   \binom{s_{d-5}+1}{2}       + \binom{d-5}{2} (2s_{d-5}-1 ) +  \binom{d-5}{3}$, 
\item [\mylabel{Y5}{\normalfont{(5)}}] $y_{d}^{d-4}=\binom{s_{d-4}}{3} +  (d-4) s_{d-4} +2\binom{d-4}{2}$,  
\item [\mylabel{Y6}{\normalfont{(6)}}] $y_{d}^{d-3}=\binom{s_{d-3}}{2}+(d-3)$,
\item [\mylabel{Y7}{\normalfont{(7)}}] $y_{d}^{d-2}=s_{d-2}$,
\item [\mylabel{Y8}{\normalfont{(8)}}] $y_{d}^{d-1}=1$,
\item [\mylabel{Y9}{\normalfont{(9)}}] $  y^{d-6}_d=\frac{1}{46080}(d^{11} - 49 d^{10} + 1060 d^9 - 13210 d^8 + 104645 d^7 - 554849 d^6\\  + 2057246 d^5 - 5611700 d^4 + 11491704 d^3 - 15859872 d^2  + 10654464 d - 967680) (d - 5)$.
\end{itemize}
The relations \ref{Y1}--\ref{Y3} are obtained by the closed formula counting partitions in dimension up to 3. Relations \ref{Y4}--\ref{Y8} are obtained by specialising  $\BL=1$  in the formulas in \cite[Prop.~7.2]{MOTIVEMOTIVEMOTIVE}. The proof of \ref{Y9} comes from combining \Cref{cor: y determinati dai primi} with our explicit computations, see the ancillary file \ancillary.
\end{example}

\subsection{Proof of \texorpdfstring{\Cref{thm: intro C}}{}}\label{sec: proof C main}
In this section we prove \Cref{thm: intro C} from the Introduction. Recall that we defined the generating function
\[
\mathsf C_x(t) =  
\begin{cases}
\displaystyle\sum_{e\geqslant \left\lceil \frac{x}{2}\right\rceil}\, {c_{e+1}^{2(e+1)-x}} t^{e-\left\lceil \frac{x}{2}\right\rceil} & \mbox{ if }x \mbox{ is even}\\ \\
\displaystyle\sum_{e\geqslant \left\lceil \frac{x}{2}\right\rceil}\, {c_{e}^{2e-x}} t^{e-\left\lceil \frac{x}{2}\right\rceil} & \mbox{ if }x\mbox{ is odd}
\end{cases}
\]  
in the ring of formal power series $\BZ\llbracket t \rrbracket$.

\begin{theorem}\label{thm:genc}
Fix $x\in\BZ_{\geqslant 0}$. Then
\begin{equation}
\label{eq:formulathm}
\Bor\left( \mathsf C_x(t) \right)=\frac{ \mathsf{c}_x(t)}{(1-2t)^{\frac{3}{2}+2x-\left\lceil\frac{x}{2}\right\rceil}}
\end{equation} 
in $\BQ\llbracket t\rrbracket$, where $\mathsf c_x(t) \in \BQ[t]$ is a polynomial of degree at most $2x-\left\lceil\frac{x}{2}\right\rceil$. In particular, $\Bor\left( \mathsf C_x(t) \right)$ is an analytic function in $t=0$.

More precisely, writing
\[
\mathsf c_x(t) = \sum_{h=0}^{2x-\left\lceil\frac{x}{2}\right\rceil} \mu_h^{(x)}t^h,
\]
we have that
\begin{itemize}
\item [\mylabel{evencase}{\normalfont{(1)}}] if $x=2\alpha$ is even, then
\[
\mu_h^{(2\alpha)} = \sum_{z=0 }^{h}\frac{(-1)^{z +h }2^h c_{z+1+\alpha}^{2z+2}}{(2z+1)!!} \sum_{y=z}^{3\alpha-1} \binom{3\alpha-y}{ h -y}\binom{y+1}{z+1} \frac{ (2y+1)!!  }{ (2y+2)!!}\left( \frac{(2y+3)(3\alpha-h)}{3\alpha-y}  -1    \right)
\]
for $0\leqslant h \leqslant 3\alpha$, and 
\item [\mylabel{oddcase}{\normalfont{(2)}}] if $x=2\alpha-1$ is odd, then 
\[
\mu_h^{(2\alpha-1)} = \sum_{z=0}^{h}\frac{(-1)^z c_{z+\alpha}^{2z+1}}{(2z+1)!!z!} \sum ^{h}_{y=z}\binom{3\alpha-2-y}{h-y}\frac{(-1)^{  y  }(2y+1)!!  2 ^{h-y}}{(y-z)!}
\]
for $0\leqslant h \leqslant 3\alpha-2$.
\end{itemize}
\end{theorem}

\begin{remark}
Before proving \Cref{thm:genc}, we comment on the statement. In  \eqref{eq:formulathm} we interpret the function 
\[
\frac{ \mathsf{c}_x(t)}{(1-2t)^{\frac{3}{2}+2x-\left\lceil\frac{x}{2}\right\rceil}}
\]
analytic at $t=0$ as its Taylor expansion at the origin.

Notice also that the statement implies a precise formula for the degree zero coefficient of $\mathsf c_x$, namely
    \[
    \mu_0^{(x)}=\begin{cases}
        c_{1+\alpha}^2&\mbox{ if }x=2\alpha,\\
        c_\alpha^1=1 &\mbox{ if }x=2\alpha-1.
    \end{cases}
    \]
\end{remark}

\begin{proofof}{\Cref{thm:genc}} 
We start proving \Cref{eq:formulathm} in the even case \ref{evencase}. If   $x=2\alpha$, our claim is equivalent to  
\[
\Bor\left(\mathsf C_{2\alpha}(t)\right) =\frac{ \mathsf{c}_{2\alpha}(t)}{(1-2t)^{\frac{3}{2}+3\alpha}}.
\] 
The case $\alpha=0$ easily follows from \Cref{lemma:borelmainprop} and \Cref{rem:limitcases}, therefore we may assume $\alpha>0$. We have
\begin{align*} 
\Bor(\mathsf C_{2\alpha}(t) )
&= \sum_{e\geqslant \alpha}\, \frac{c_{e+1}^{2(e+1-\alpha)}}{(e-\alpha)!}t^{e-\alpha} \\
&= \sum_{e\geqslant \alpha}\, \frac{t^{e-\alpha} }{(e-\alpha)!} \sum_{y=\alpha }^{\min(e+1,4\alpha)}\frac{(2(e+1-\alpha))!}{(2(e+1-y))!!(2(y-\alpha))!}d_{y}^{2(y-\alpha)}\\
& = \sum_{y=\alpha+1 }^{4\alpha} \sum_{e\geqslant y-1}\, \frac{t^{e-\alpha} }{(e-\alpha)!}\frac{(2(e+1-\alpha))!}{(2(e+1-y))!!(2(y-\alpha))!}d_{y}^{2(y-\alpha)}   \\ 
& = \sum_{y=1 }^{3\alpha} \sum_{e\geqslant y}\, \frac{t^{e-1} }{(e-1)!} \frac{(2e)!}{(2(e-y))!!(2y)!}d_{y+\alpha}^{2y}   \\ 
& = \sum_{y=1 }^{3\alpha} \sum_{e\geqslant y}\sum_{z=0}^{y}\, \frac{t^{e-1} }{(e-1)!}\frac{(-1)^z(2e)!}{(2e-2y)!!(2z)!!(2y-2z)!}c_{y-z+\alpha}^{2y-2z},
\end{align*}
where we have used \Cref{prop: inversion}\ref{inversion-3} in the last step. Noting that the term $z=0$ does not contribute and performing a shift of the variables, we continue the above chain of identities to get
\begin{align*} 
\Bor(\mathsf C_{2\alpha}(t) )
&= \sum_{y=1 }^{3\alpha} \sum_{e\geqslant y}\sum_{z=1}^{y}\,  \frac{(-1)^{z+y}e(2e-1)!!}{(e-y)!(y-z)!(2z-1)!!z!}c_{z+\alpha}^{2z} t^{e-1}   \\ 
& = \sum_{y=1 }^{3\alpha} \sum_{z=1}^{y}\frac{(-1)^{z+y}c_{z+\alpha}^{2z}t^{y-1}}{(y-z)!(2z-1)!!z!}\sum_{e\geqslant 0}\,  \frac{(e+y)(2e+2y-1)!!}{e!} t^{e}   \\ 
& = \sum_{y=1 }^{3\alpha} \sum_{z=1}^{y}\frac{(-1)^{z+y}c_{z+\alpha}^{2z}t^{y-1}}{(y-z)!(2z-1)!!z!}\frac{(2y-1)!!(y+t)}{(1-2t)^{\frac{3}{2}+y}},
\end{align*}
where we have used \Cref{lemma:techborel} for the last identity. Let us ignore the factor $(1-2t)^{-\frac{3}{2}-y}$ from now on. What remains can be rewritten as 
\begin{align*}
\sum_{z=1 }^{3\alpha} \frac{(-1)^{z }c_{z+\alpha}^{2z}}{(2z-1)!!z!}&\sum_{y=z}^{3\alpha}\frac{(-1)^{ y}(2y-1)!!t^{y-1}(y+t)(1-2t)^{3\alpha-y}}{(y-z)!}  \\  
& =\sum_{z=1 }^{3\alpha} \frac{(-1)^{z }c_{z+\alpha}^{2z}}{(2z-1)!!z!}\sum_{y=z}^{3\alpha}\sum_{h=0}^{3\alpha-y}\binom{3\alpha-y}{h}\frac{(-1)^{ 3\alpha+h}(2y-1)!!t^{3\alpha-h-1}(y+t)2^{3\alpha-y-h}}{(y-z)!}  \\ 
& = \sum_{z=0 }^{3\alpha-1} \frac{(-1)^{z }c_{z+1+\alpha}^{2z+2}}{(2z+1)!!(z+1)!}\sum_{y=z}^{3\alpha-1}\sum_{h=y}^{3\alpha-1}\binom{3\alpha-y-1}{ h -y}\frac{(-1)^{ h}(2y+1)!!(y+1+t)2^{h-y}}{(y-z)!}t^{h}.
\end{align*}
Now, for the purpose of extracting the coefficient of $t^h$, we split the $(y+1+t)$ term appearing in the last fraction. This gives
\begin{multline*}
\sum_{z=0 }^{3\alpha-1} \frac{(-1)^{z }c_{z+1+\alpha}^{2z+2}}{(2z+1)!!(z+1)!}\sum_{y=z}^{3\alpha-1}\Biggl(\sum_{h=y}^{3\alpha-1}\binom{3\alpha-y-1}{ h -y}\frac{(-1)^{ h}(2y+1)!!(y+1 )2^{h-y}}{(y-z)!}t^{h}  \\  
+\sum_{h=y}^{3\alpha-1}\binom{3\alpha-y-1}{ h -y}\frac{(-1)^{ h}(2y+1)!! 2^{h-y}}{(y-z)!}t^{h+1}\Biggr).
\end{multline*}
Rearranging, the last expression equals
\begin{align*}
&\sum_{z=0 }^{3\alpha-1} \sum_{y=z}^{3\alpha-1} \sum_{h=y}^{3\alpha }\frac{(-1)^{z }c_{z+1+\alpha}^{2z+2}}{(2z+1)!!(z+1)!}\binom{3\alpha-y}{ h -y}\frac{(-1)^{ h}(2(y+1)(3\alpha-h)-h+y)(2y+1)!! 2^{h-y-1}}{(3\alpha-y)(y-z)!}t^{h} \\ 
&= \sum_{h=0}^{3\alpha } \sum_{z=0 }^{h}\sum_{y=z}^{3\alpha-1} \frac{(-1)^{z }c_{z+1+\alpha}^{2z+2}}{(2z+1)!!(z+1)!} \binom{3\alpha-y}{ h -y}\frac{(-1)^{ h}(2(y+1)(3\alpha-h)-h+y)(2y+1)!! 2^{h-y-1}}{(3\alpha-y)(y-z)!}t^{h}.
\end{align*}
The coefficient $\mu_h^{(2\alpha)}$ of $t^h$ in the above sum can now be rewritten as  
\begin{align*}
    \mu_h^{(2\alpha)}&=\sum_{z=0 }^{h}\frac{(-1)^{z +h }2^h c_{z+1+\alpha}^{2z+2}}{(2z+1)!!} \sum_{y=z}^{3\alpha-1} \binom{3\alpha-y}{ h -y}\binom{y+1}{z+1} \frac{ (2y+1)!!  }{ (2y+2)!!}\left( \frac{(2y+3)(3\alpha-h)}{3\alpha-y}  -1    \right).
\end{align*} 
This finishes the proof of the even case \ref{evencase}.

\smallbreak
It remains to prove \Cref{eq:formulathm} in the odd case \ref{oddcase}.
If  $x=2\alpha-1$, our claim is equivalent to
\[
\Bor\left(\mathsf C_{2\alpha-1}(t)\right) =\frac{ \mathsf{c}_{2\alpha-1}(t)}{(1-2t)^{3\alpha-\frac{1}{2}}} .
\] 
We compute
\begin{align*}
\Bor(\mathsf C_{2\alpha-1}(t)) 
&=\sum_{e\geqslant \alpha}\, \frac{c_{e}^{2e-2\alpha+1}}{(e-\alpha)!} t^{e-\alpha}  \\
&= \sum_{e\geqslant \alpha}\, \frac{t^{e-\alpha} }{(e-\alpha)!} \sum_{y=\alpha }^{\min(e,4\alpha-2)}\frac{(2e-2\alpha+1)!}{(2e-2y )!!(2y-2\alpha+1)!}d_{y}^{2y-2\alpha+1} \\
&=\sum_{y=\alpha }^{ 4\alpha-2 } \sum_{e\geqslant y}\, \frac{t^{e-\alpha} }{(e-\alpha)!} \frac{(2e-2\alpha+1)!}{(2e-2y )!!(2y-2\alpha+1)!}d_{y}^{2y-2\alpha+1} \\
&=\sum_{y=0 }^{ 3\alpha-2 } \sum_{e\geqslant y}\sum ^{y}_{z=0}(-1)^{z+y}\, \frac{t^{e} }{e!} \frac{(2e+1)!}{(2e-2y )!!(2y-2z)!!(2z+1)!}c_{z+\alpha}^{2z+1}\\
&=\sum_{y=0 }^{ 3\alpha-2 } \sum ^{y}_{z=0}\frac{(-1)^{z+y}(2y+1)!!t^yc_{z+\alpha}^{2z+1}}{(y-z)!(2z+1)!!z!}\sum_{e\geqslant 0}\,  \frac{(2e+2(y+1)-1)!!}{(2(y+1)-1)!!e!}t^{e} \\
&=\sum_{y=0 }^{ 3\alpha-2 } \sum ^{y}_{z=0}\frac{(-1)^{z+y}(2y+1)!!t^yc_{z+\alpha}^{2z+1}}{(y-z)!(2z+1)!!z!}\frac{1}{(1-2t)^{\frac{3}{2}+y}}\\
&=\frac{1}{(1-2t)^{3\alpha-\frac{1}{2}}}\sum_{y=0 }^{ 3\alpha-2 } \sum ^{y}_{z=0}\sum_{h=0}^{3\alpha-2-y}\binom{3\alpha-2-y}{h}\frac{(-1)^{z+y}(2y+1)!!t^y(-2t)^{3\alpha-2-y-h}c_{z+\alpha}^{2z+1}}{(y-z)!(2z+1)!!z!}\\
&=\frac{1}{(1-2t)^{3\alpha-\frac{1}{2}}}\sum_{h=0 }^{ 3\alpha-2 }\sum_{y=0}^{h} \sum ^{y}_{z=0}\binom{3\alpha-2-y}{h-y}\frac{(-1)^{  h +z}(2y+1)!!(-2)^{h-y}c_{z+\alpha}^{2z+1}}{(y-z)!(2z+1)!!z!}t^{ h }
\end{align*} 
The coefficient $\mu_h^{(2\alpha-1)}$ of $t^h$ in the above sum can now be rewritten as  
\begin{align}
\mu_h^{(2\alpha-1)}&=\sum_{z=0}^{h}\frac{(-1)^z c_{z+\alpha}^{2z+1}}{(2z+1)!!z!} \sum ^{h}_{y=z}\binom{3\alpha-2-y}{h-y}\frac{(-1)^{y}(2y+1)!!  2 ^{h-y}}{(y-z)!} .
\end{align}
This completes the proof of the odd case \ref{oddcase}, and hence of the theorem.
\end{proofof} 

\begin{remark}
\label{rmk:hypergeom}
For each $h = 0,\ldots,3\alpha-2$ there is an identity\footnote{Obtained with the help of Mathematica \cite{Mathematica}.} 
\[
\mu_h^{(2\alpha-1)} = \sum_{z=0}^{h} \binom{3\alpha-z-2}{h-z}\frac{2^{h-z}}{z!}{_2\mathrm{F}_1}\left(\frac{3}{2}+z,z-h;2+z-3\alpha;-1\right) c_{z+\alpha}^{2z+1},
\]   
where ${_2\mathrm{F}_1}(a,b;c;t)$ is the ordinary hypergeometric function, see e.g.~\cite{OD_Hypergeom}.
\end{remark}
\section{M-partitions: reduction to socle in degree \texorpdfstring{$\geqslant 3$}{}}\label{sec: algor}
In this section we go one step further in our socle reduction technique: we provide a recursion (cf.~\Cref{prop:recur}) that recovers the numbers $c_\bullet^\bullet$ from the number of partitions with socle in degree at least 3. The key character in this recursion is the number of \emph{M-partitions}, that we introduce in \Cref{def:Mpartition} below. In \Cref{subsec:algorithms}  we discuss explicit methods to effectively compute these numbers.

\begin{definition}
\label{def:Mpartition}
For a triple of positive integers $k, q, m>0$, we call \emph{M-partition} of type $(k,q,m)$ a $(k-1)$-dimensional partition $\lambda\in \mathrm{P}^k_{1+k+q+m}$ such that 
\begin{itemize}
    \item [(i)] $\Soc(\lambda)\subset \lambda_{\geqslant 3}$,
    \item [(ii)]  $h_\lambda(1)=k$,
    \item [(iii)]  $h_\lambda(2)=q$,
    \item [(iv)]  $\sum_{i\geqslant 3} h_\lambda(i)=m$.
\end{itemize}
We denote by $\alpha_{q,m}^k$ the number of M-partitions of type $(k,q,m)$. We formally set $\alpha_{0,0}^0=1$, and $\alpha_{q,0}^k=0$ for $k,q>0$.
\end{definition}
The next result shows that the number of quadrics in an M-partition is always at least the number of variables.
\begin{prop}\label{prop:minq}
We have $\alpha_{q,m}^k =0$, for $q < k$.
\end{prop}

\begin{proof}
Let $\lambda$ be an M-partition of type $(k,q,m)$ and let $I$ be the monomial ideal associated to $\lambda$. Identifying lattice points in $\lambda$ with monomials, we set
\[
\lambda_{=1}=\{x_1, \dots, x_k\},
\]
and 
\[
\lambda_{=2}=\coprod_{i=1}^4 Q_i
\]
to be the decomposition of degree 2 elements given by (up to reordering of the variables)
\begin{align*}
    Q_1&=\set{x_1^2, \dots, x_l^2}, \quad \mbox{for some } 1\leqslant l\leqslant k,\\
    Q_2&= \set{x_a x_b\in \lambda| a,b=1, \dots, l \mbox{ and } a\neq b},\\
    Q_3&= \set{x_a x_b\in \lambda| a=1, \dots, l \mbox{ and } b=l+1, \dots, k},\\
    Q_4&= \set{x_a x_b\in \lambda| a,b=l+1, \dots, k \mbox{ and } a\neq b}.
\end{align*}
Since $\lambda$ admits no socle in degree 1, all the variables in $\lambda_{=1}$ are obtained as derivates of some monomials in $\lambda_{=2}$. In particular, $Q_1$ contributes precisely with $l$ variables, $Q_2$ does not contribute with new variables (with respect to the ones already introduced by $Q_1$), and $Q_3$ contributes with at most $|Q_3|$ new variables. 

On the other hand, $\lambda$ admits no socle in degree 2, which implies that for all monomials $x_ax_b\in Q_4$, there must exist a third different variable $x_c $ such that $x_ax_c, x_bx_c\in Q_4 $. This readily implies that $Q_4$ contributes with at most $|Q_4|$ new variables. 

Summing up all the contributions yields $q\geqslant k$, which implies the result.
\end{proof}

\begin{remark}
\label{rem:alphasprop}
The numbers $\alpha_{q,m}^k$ may be refined in several ways. For instance, set
\[
\alpha_{q,m,\ell}^k=\left|\Set{\lambda\in \mathrm{P}^k_{1+k+q+m} | \lambda \mbox{ is an M-partition of type } (k,q,m) \mbox{ and } \ell(\lambda)=\ell }\right|.
\]
The numbers $\alpha_{q,m,\ell}^k$ play a central role in the counting procedure that we develop in \Cref{subsec:algorithms}.
For a more refined example, let $h$ be a Hilbert--Samuel function. Define
\[
\alpha_{h}=\left|\Set{\lambda\in \mathrm{P}^{h(1)}_{\lvert {h}\rvert} | \lambda \mbox{ is an M-partition of type } \left(h(1),h(2),\sum_{i\geqslant 3} h(i)\right) \mbox{ and } h_\lambda=h }\right|.
\]
We  compute explicitly the numbers $\alpha_{q,m,\ell}^k$ and $\alpha_{h}$ in \Cref{prop:max1,prop:max2,prop:max3,thm:expexp}, where we study some families of compressed and anti-compressed partitions, see \Cref{subsec:compressed}. 
\end{remark}

\begin{example}
The number of M-partitions of type $(3,5,13)$ and  length $4$ is $\alpha_{5,13,4}^3=531$.  Precisely, this number decomposes as
\[
\alpha_{5,13,4}^3=\alpha_{h_1}+\alpha_{h_2}=504+27,
\]
where  $h_1=(1, 3, 5, 7, 6)$ and $h_2=(1, 3, 5, 6, 7)$. On the other hand, we have in total $\alpha_{5,13}^3 = 43260$.
\end{example}

\begin{definition}
\label{def: c molti indici}
Given a triple of positive integers $k,q,m$, and a pair of nonnegative integers $e,a\geqslant 0$, we denote by
$ c_{e,(k,q,m)}^a$ the number of partitions $\lambda\in \mathrm{C}_e^a$ such that the partition
\[
\left( \Soc(\lambda)\cap \lambda_{\geqslant 3} \right)^\perp
\]
is an M-partition of type $(k,q,m)$. For $a>0$, we set $c_{e, (0,0,0)}^a$ to be the number of partitions $\lambda\in\mathrm{C}_e^a $ such that $\Soc(\lambda)= \lambda_{=2}$.
\end{definition}

In the next result, we show that  knowing the numbers  $\alpha_{\bullet, \bullet}^\bullet$  implies the knowledge of the numbers $c_{\bullet}^\bullet$.

\begin{prop}
\label{prop:recur}
The numbers $c_{e,(k,q,m)}^a$ satisfy the relations

\begin{equation}\label{eq:recurc}
\begin{aligned}
c_{e,(k,q,m)}^0 &= 0, \\
c_{e,(k,q,m)}^a   &= \alpha_{q,m}^k\binom{a}{k} \binom{ \binom{a+1}{2}-q}{e-q-m} - \sum_{i=1}^{a-1}{\binom{a}{i}c_{e,(k,q,m)}^i}. 
\end{aligned}
\end{equation}
In particular, there is an identity
\[
c_e^a = \sum_{\substack{k \leqslant a \\ q+m\leqslant  e}} c_{e,(k,q,m)}^a.
\]
\end{prop}

\begin{proof}
The vanishing $c_{e,(k,q,m)}^0 = 0$ is a direct consequence of the definition of the numbers $c_{e,(k,q,m)}^a$. To construct a partition $\lambda$ among the ones counted by $c_{e,(k,q,m)}^a$, we first choose $k$ variables among the $a$ variables given, and then choose an M-partition $\pi$ of type $(k,q,m)$ in the $k$ variables we selected. This amount to a contribution of 
\[
\alpha_{q,m}^k\binom{a}{k}.
\]
Next, we have to add $ e-q-m$ boxes in degree 2 to $\pi$, which contributes 
    \[
    \binom{ \binom{a+1}{2}-q}{e-q-m}.
    \]
Notice that not all partitions constructed in this way have embedding dimension $a$. In fact, they could be represented by some partitions counted by $ c_{e,(k,q,m)}^i$, for some $i<a$. Summing it up yields
  \begin{align*}
c_{e,(k,q,m)}^a   &= \alpha_{q,m}^k\binom{a}{k} \binom{ \binom{a+1}{2}-q}{e-q-m} - \sum_{i=1}^{a-1}{ \binom{a}{i}c_{e,(k,q,m)}^i  }. 
\end{align*}
The second claim follows by subdiving the set $\mathrm{C}_e^k$ according to the type $(k,q,m)$ of the M-partition $ \left( \Soc(\lambda)\cap \lambda_{\geqslant 3} \right)^\perp$. 
\end{proof}

The recursion in \Cref{prop:recur} shows that only finitely many $c_{e,(k,q,m)}^a$ contribute to the computation of $c_{e}^a$. It is also worth mentioning that some of the contibuting  $c_{e,(k,q,m)}^a$ might be zero: for instance, we have $c_{3,(2, 2, 1)}^3 = 0$.

\begin{example}
Figures \ref{fig:c(e,3mq,a)} and \ref{fig:c(e,2mq,a)} depict two examples of partitions  $\lambda\in \mathrm{P}^3_{16}$ and $\lambda'\in \mathrm{P}^3_{13}$ counted by $c_{12,(3,4,8)}^3$ and $c_{9,(2,3,5)}^3$ respectively. 
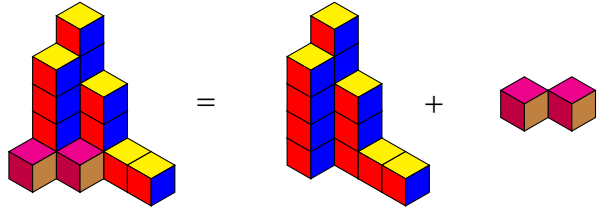
\begin{figure}[ht]
    \centering
    \begin{tikzpicture}
           
\node at (0,0) { \begin{tikzpicture}[scale=0.36] 
\planepartitionone{{5,3,1,1},{4}} 
\end{tikzpicture}} ;
\node at (1.5,0) {=}; 
\node at (6,0){\begin{tikzpicture}[scale=0.36] 
\planepartitionone{} 
\end{tikzpicture}} ;
\node at (4.5,0) {+}; 
\node at (3.5,0 ){ \begin{tikzpicture}[scale=0.36] 
\planepartition{{5,3,1,1},{4}}    
\end{tikzpicture}}; 
\end{tikzpicture} 
    \caption{A partition $\lambda\in \mathrm{P}^3_{16}$ counted by $c_{12,(3,4,8)}^3$.}
    \label{fig:c(e,3mq,a)}
\end{figure} 
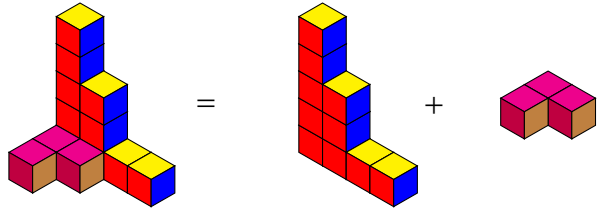
\begin{figure}[ht]
    \centering
    \begin{tikzpicture}
           
\node at (0,0) { \begin{tikzpicture}[scale=0.36] 
\planepartitiondue{{5,3,1,1}} 
\end{tikzpicture}} ;
\node at (1.5,0) {=}; 
\node at (6,0){\begin{tikzpicture}[scale=0.36] 
\planepartitiondue{} 
\end{tikzpicture}} ;
\node at (4.5,0) {+}; 
\node at ( 3.5,0 ){ \begin{tikzpicture}[scale=0.36] 
\planepartition{{5,3,1,1}}    
\end{tikzpicture}}; 
\end{tikzpicture} 
    \caption{A partition $\lambda'\in \mathrm{P}^3_{13}$ counted by $c_{9,(2,3,5)}^3$.}
    \label{fig:c(e,2mq,a)}
\end{figure}
\end{example}

\section{Exact count of M-partitions} \label{sec:exactformulas}
In this section, we provide exact formulas for the generating functions counting partitions with prescribed Hilbert--Samuel functions and constraints on their socles. More precisely, we find generating functions for some collections of numbers $\alpha^\bullet_{\bullet, \bullet, \bullet}$ (cf.~\Cref{rem:alphasprop}).
Finding such generating functions reduces the complexity of the computations of the numbers $c_\bullet^\bullet$ (and, in turn, of the number of partitions $p_\bullet^\bullet$, via \Cref{prop: inversion}) and offers interesting perspectives on the asymptotic distribution of partitions.

\subsection{Compressed and anti-compressed partitions}
\label{subsec:compressed}
Fix $n,r\in \BZ_{\geqslant 1}$ and a vector $e \in \BN^r$.   We say that an $(n-1)$-dimensional partition $\lambda$ is \emph{compressed} (resp.~\emph{anti-compressed}) with respect to $e$ if it has maximal (resp.~minimal) size among the partitions with    socle type $e_\lambda=e$ and embedding dimension $h_\lambda(1)=n$.

\begin{prop}
\label{prop:max1}
Fix $n\in \BZ_{\geqslant 1}$. A partition $\lambda$ with embedding dimension $h_\lambda(1)=3n$ and socle type $(0,0,0,n)$ has  Hilbert--Samuel function $(1,3n,3n,n)$. In particular, it is a compressed M-partition of type $(3n,3n,n)$ and length $3$. Moreover, the number $\alpha_{3n,n,3}^{3n}$ of such partitions is  
\begin{equation}
\label{eqn:alpha-3n}
\alpha_{3n,n,3}^{3n}=\frac{(3 n)!}{6^nn!}.
\end{equation}
\end{prop}

\begin{proof}
If $n$ cubics generate $3n$ linear terms, it means that  the set of cubics in $\lambda$ consists of   $x_1y_{1}z_1,\dots, x_n y_n z_n$, which in particular means that there are $3n$ quadrics, all of the form $x_iy_i, x_iz_i, y_iz_i$, and that the associated partition is compressed. Therefore, the possible partitions of this form correspond to the numbers of partitions of a set of $3n$ elements in $n$ subsets of 3 elements, that is
\[ 
\frac{1}{n!}\binom{3n}{3}\binom{3n-3}{3}\cdots \binom{6}{3} \binom{3}{3}=\frac{(3 n)!}{6^nn!}.\qedhere
\] 
\end{proof}

\begin{prop} 
\label{prop:max2} 
Fix $n\in \BZ_{\geqslant 1}$. A partition $\lambda$ with embedding dimension $h_\lambda(1)=3n-1$ and socle type $(0,0,0,n)$ has  possible Hilbert--Samuel functions
\begin{itemize}
    \item $(1,3n-1,3n,n)$,
    \item $(1,3n-1,3n-1,n)$.
\end{itemize}
In particular, in the first case it is a compressed M-partition of type $(3n-1,3n,n)$ and length 3, and the number  
of such partitions is  
\[
\alpha_{3n,n,3}^{3n-1}=\frac{3(n-1)}{2}\alpha_{3n,n,3}^{3n}. 
\] 
In the second case $\lambda$ is anti-compressed, and the number of such partitions is 
\[
\alpha_{3n-1,n,3}^{3n-1}=2\alpha_{3n,n,3}^{3n}. 
\]
\end{prop}

\begin{proof}
If $n$ cubics generate $3n-1$ linear terms, the cubics are (up to symmetry or relabelling of the variables) of one of the two forms
\begin{itemize}
    \item $x_1y_{1}z_1,\dots,x_{n-1}y_{n-1}z_{n-1}, x_n y_n x_i$, for some $i=1, \dots, n-1$,
    \item $x_1y_{1}z_1,\dots,x_{n-1}y_{n-1}z_{n-1}, x_n y^2_n$.
\end{itemize}
In the first case, the $n$ cubics generate $3n$ quadrics and the partition is compressed, while in the second case they generate $3n-1$ quadrics and the partition is anti-compressed. A simple combinatorial computation yields that in the compressed case the number of such partitions is
\[
\frac{1}{2}\binom{3n-1}{2}\cdot(3n-3)\cdot\frac{(3n-3)!}{6^{n-1}(n-1)!}=\frac{3(n-1)(3n-1)!}{4\cdot 6^{n-1}(n-1)!},
\]
while in the anti-compressed case is 
\[
(3n-1)(3n-2)\cdot\frac{(3n-3)!}{6^{n-1}(n-1)!}=\frac{(3n-1)!}{ 6^{n-1}(n-1)!}.
\]
The result now follows from \Cref{eqn:alpha-3n}.
\end{proof}

\begin{prop}
\label{prop:max3} 
Fix $n\in \BZ_{\geqslant 1}$. An anti-compressed partition $\lambda$ with embedding dimension $h_\lambda(1)=3n-2$ and socle type $(0,0,0,n)$  has  Hilbert--Samuel function $(1,3n-2,3n-2,n)$.
The number of such partitions is
\[
\alpha_{3n-2,n,3}^{3n-2}=(3n-2)^2\alpha_{3(n-1),n-1,3}^{3(n-1)}.
\] 
\end{prop}

\begin{proof}
A case by case analysis shows that the possible collections of cubics in $\lambda$ (up to symmetry or relabelling of the variables) are of one of the two forms
\begin{itemize}
    \item $x_1y_{1}z_1,\dots,x_{n-1}y_{n-1}z_{n-1}, x^3_n$,
    \item $x_1y_{1}z_1,\dots,x_{n-2}y_{n-2}z_{n-2}, x_{n-1}y^2_{n-1}, x_{n}y^2_{n}$.
\end{itemize}
   The first configuration contributes
    \[
    (3n-2)\alpha_{3(n-1),n-1,3}^{3(n-1)}.
    \]
    The second configuration contributes
    \[
    12\binom{3n-2}{4}\alpha_{3(n-2),n-2,3}^{3(n-2)}=12\binom{3n-2}{4}\frac{(3n-6)!}{6^{n-2}(n-2)!}=3 \frac{(3n-2)!}{6^{n-1} \cdot(n-2)!}.
    \]
    Summing the contributions yields the result.
\end{proof}

We provide a closed formula for the generating functions counting infinite families of anti-compressed partitions of arbitrary length. These partitions are a special instance of the \emph{hydral} partitions, cf.~\cref{def: hydral}.
\begin{theorem}\label{thm:expexp}
Fix a positive integer $s\geqslant 3$ and define the Hilbert--Samuel function 
\[
h_{n,s}=(1,\underbrace{n,n,\ldots ,n}_{s}).
\]
Then we have the identity
\[
\alpha_{h_{n,s}}=\sum_{k=0}^n\binom{n}{k}k^{n-k}.
\]
In particular, $\alpha_{h_{n,s}}$ is independent of $s$ and  there is an identity
\[
\Bor\left( \sum_{n\geqslant 0} \alpha_{h_{n,s}} t^n\right) = e^{te^t}.
\]
\end{theorem}

\begin{proof}
The second claim easily follows from the first one, by noticing that 
\[
\Coef_{t^n}{\left(\IB\left(  e^{te^t}\right)\right)}=\sum_{k=0}^n\binom{n}{k}k^{n-k},
\]
which is proved by a direct computation. We divide the proof of the first claim in several steps.

\underline{Step I.}  For $s\geqslant 3$,  set
\[
\mathrm{P}_{n,s}=\Set{\lambda\in \mathrm{P}_{ns+1}^n | \Soc(\lambda)\subset \lambda_{\geqslant 3} \mbox{ and } h_\lambda=h_{n,s}}
\]
to be the collections of partitions of interest counted by $\alpha_{h_{n,s}}$ and 
\[
\mathrm{P}_{n}^{\vdash}=\Set{(\Lambda, \overline{a}) |  \Lambda\subset 2^{[n]}, \,  \coprod_{A\in \Lambda} A=[n], \,   \overline{a}=(a_A)_{A\in \Lambda}\in\prod_{A\in\Lambda} A } 
\]
to be the collection of partitions of the set  $[n]=\set{1,\ldots,n}$ with a \emph{coloured} element in each subset of the partition. Consider the map
\begin{equation}\label{eqn: forma delle partizioni step 3}
    \begin{tikzcd}[row sep=tiny] \mathrm{P}_{n}^{\vdash}\arrow[r,"\Psi_s"]&\mathrm{P}_{n,s}\\
(\Lambda, \overline{a})\arrow[r,mapsto]& \left(\bigcup_{A\in \Lambda} \Set{x_{a_A}^{s-1}x_b  |  b\in A}\right)^\perp
\end{tikzcd}
\end{equation}
where we denote  partitions in $\mathrm{P}_{n,s}$ using monomial ideals through the correspondence of \Cref{sec: monomial ideals}. The map $\Psi_s $ is well-defined, as $h_{\Psi_s(\Lambda, \overline{a})}=h_{n,s}$. Notice that the number of partitions of $[n] $ with a chosen coloured element in each subset of the partition is given by 
\[
\lvert \mathrm{P}_{n}^{\vdash}\rvert = \sum_{k=0}^n\binom{n}{k}k^{n-k}.
\]
Therefore, the conclusion of the theorem folllows by proving that $\Psi_s$ is a bijection.

\underline{Step II.} We prove that  $ \Psi_{s-1}$ is a bijection by induction on $s\geqslant 3$. In this step, we prove the inductive step. Assume therefore that $\Psi_s$ is a bijection. Consider the map
\[
\begin{tikzcd}[row sep=tiny]
\mathrm{P}_{n,s}\arrow[r,"\tau_s"]&\mathrm{P}_{n,s-1}\\
    I\arrow[r,mapsto]&I+\mathfrak m^s,
\end{tikzcd}
\]
    where $\mathfrak m\subset\BC[x_1, \dots, x_n] $ is the monomial
    maximal ideal. We claim that $\tau_s $ is a bijection. 

    For surjectivity, by the inductive assumption, we have that all ideals $J\in \mathrm{P}_{n,s-1}$ are of the form  
    \begin{equation}
    \label{eq:ideal J}
         J=\left(\bigcup_{A\in \Lambda} \Set{x_{a_A}^{s-2}x_b  |  b\in A}\right)^\perp
    \end{equation}
    for some $(\Lambda, \overline{a})\in \mathrm{P}_{n}^{\vdash}$. However, the ideal $J$ satisfies
    \[
    J=\Psi_s\left(\left(\bigcup_{A\in \Lambda} \Set{x_{a_A}^{s-1}x_a  |  a\in A} \right)^\perp \right),
    \]
    which yields surjectivity. 

    For injectivity, assume there are two ideals $I, I'\in \mathrm{P}_{n,s}$ such that 
    \[
    I+\mathfrak m^s=I'+\mathfrak m^s.
    \]
Since by induction we know that all ideals in $ \mathrm{P}_{n,s-1}$ are of the form \eqref{eq:ideal J}, without loss of generality we can suppose
\[
I=\left(\bigcup_{A\in \Lambda} \Set{x_{a_A}^{s-1}x_a  |  a\in A}\right)^\perp,
\]
for some $(\Lambda, \overline{a})\in \mathrm{P}_{n}^{\vdash}$. Assume by contradiction that $I\neq I'$. Then it means that $\BC[x_1, \dots, x_n]/I'$ contains a monomial of degree $s$ which is not of the form $x^{s-1}_{a_A}x_b$ for some $A\in \Lambda$ and $b\in A$. The only possible choice would be a monomial of the form $x^{s-2}_{a_A}x^2_b$, which is a contradiction since it would imply that $I'+\mathfrak m^s $ contains $x^{s-3}_{a_A}x^2_b$ as well. 

\underline{Step III.} Now we prove the base of the induction. Let $I\in  \mathrm{P}_{n,3}$ be an ideal corresponding to a partition with Hilbert--Samuel function $(1,n,n,n)$ and with socle in degree 3. Set $R=\BC[x_1, \dots, x_n]$. Surjectivity of $\Psi_3$ is equivalent to the claim that 
    \begin{align}\label{eqn: monomials form proof}
         R/I\cong \bigoplus_{A\in \Lambda}\bigoplus_{b\in A}\BC\cdot x^2_{a_A}x_b,
    \end{align}
for some $(\Lambda, \overline{a})\in  \mathrm{P}^{\vdash}_{n}$, which we now prove. Up to relabelling the variables, let 
    \[
    C_0=\set{x_1^3, \dots, x_k^3}\subset R/I
    \]
    be the third-power cubic monomials in $R/I$, for some $k\leqslant n$. This implies that $R/I$ contains the quadratic terms 
    \[Q_0=\set{x_1^2, \dots, x_k^2}\subset R/I\]
    and the linear terms
    \[L_0=\set{x_1, \dots, x_k}\subset R/I.\]
The collection of these monomials (and the constant term 1) contributes $(1,k,k,k)$ to the Hilbert--Samuel function. To realise  a possible partition of the required form, we still have to add $n-k$ different cubic monomials $m_1, \dots, m_{n-k}$ to $R/I$. Set 
\[
C_i=C_{i-1}\cup\set{m_i}, \quad i=1, \dots, n-k.
\]
Adding a cubic monomial $m_i$ to $R/I$ implies that at most other 3 quadratic monomials $q_{i, 1},q_{i,2}, q_{i, 3}$ and linear terms $y_{i, 1}, y_{i,2}, y_{i,3}$ have to be in  $R/I$ as well. Set
\begin{align*}
    Q_i&=Q_{i-1}\cup\set{q_{i,1}, q_{i,2}, q_{i,3}}, \quad i=1, \dots, n-k,\\
    L_i&=L_{i-1}\cup\set{y_{i,1}, y_{i,2}, y_{i,3}}, \quad i=1, \dots, n-k,
\end{align*}
Notice that at the `step' 0, the number of monomials in degree $0,1,2,3$ (i.e.~the cardinalities of $\set{1}, L_0, Q_0, C_0$) is
    \[
    (1,k,k,k),
    \]
and we ask that at the `step' $n-k$ the number of monomials in degree  $0,1,2,3$ (i.e.~the cardinalities of $\set{1}, L_{n-k}, Q_{n-k}, C_{n-k}$) to be
\[
(1,n,n,n).
\]
Now, at each `step' $i$, adding a cubic $m_i$ implies adding some quadratic and linear monomial, therefore contributes with one of the possible following vectors
     \begin{align*}
        e_1&=(0,0,0,1)\\
        e_2&=(0,0,1,1)\\
        e_3&=(0,1,1,1)\\
        e_4&=(0,0,2,1)\\
        e_5&=(0,1,2,1)\\
        e_6&=(0,2,2,1)\\
        e_7&=(0,0,3,1)\\
         e_8&=(0,2,3,1)\\
        e_9&=(0,3,3,1)
    \end{align*}
    to the cardinalities of $\set{1}, L_{i}, Q_{i}, C_{i}$.  This means that performing all the `steps' for $i=1, \dots, n-k$ yields a relation
\begin{align*}
    (1,k,k,k) + \sum_{i=1}^{n-k}e_{j_i}  =(1,n,n,n).
\end{align*}
Since each contribution vector $e_i$ could appear with multiplicity, the relation can be rewritten as 
\begin{align*}
    (1,k,k,k) + \sum_{i=1}^{9} a_i e_{i}  =(1,n,n,n),
\end{align*}
    for some nonnegative integers $a_i\in \BZ_{\geqslant 0}$. Solving the linear system, we find that 
    \begin{align*}
      a_2&=a_4=a_5=a_7=a_8=0,\\
      a_1&=a_6+2a_9\\
      a_3&=n-k-2a_6-3a_9.
    \end{align*} 
    This in particular means that no `step' could contribute with vectors $ e_2, e_4, e_5, e_7, e_8$. We claim now that $a_  1=0$. This would conclude the proof, since it would imply that each step contributes with $e_3=(0,1,1,1)$, which is easily seen to be equivalent to adding at step $i$ a cubic monomial of the form $x_a^2 y$, where $a=1, \dots, k$ and $y$ is a new variable. This in fact contributes the new monomials $x_ay, y$ and would imply that our partition is of the form \eqref{eqn: monomials form proof}.

    To prove that $a_1=0$, we have to prove that no step could contribute with $e_1$. Therefore, suppose the step $i$ is the first one to contribute with $e_1$.  There are two cases: $m_i$  could either be of the form $xyz$ for three different variables $x,y,z$ or  of the form $x^2y$, for two different variables $x,y$. We spell out the proof only that the first case; the second case follows by an analogous reasoning.

    So say that $m_i=xyz$ contributes with $(0,0,0,1)$ at step $i$. This means that $xy,yz,xz\in Q_{i-1} $ and $x,y,z\in L_{i-1}$ were already contained in the step $i-1$. Since clearly $xy,yz,xz\notin Q_{0}$, there should have been three intermediate steps $1\leqslant i'< i''< i'''<i$ adding cubic monomials of the form 
    \begin{align*}
        m_{i'}&=xyt,\\
         m_{i''}&=xzu,\\
          m_{i'''}&=yzw,
    \end{align*}
for some $t,u,w$. We require these steps to be respectively the first time that the quadric monomials $xy, xz, yz$ appear at the corresponding step. In particular, the step $i'''$ contributes with quadrics $yz, yw, zw $ and variables $y,z,w$. Since $y,z\in L_{i'''-1}$, at step $i'''$ we could add either 0 or 1 new linear term, so the only possibility is that this step contributes with $e_4$ and therefore $w\notin L_{i'''-1}$. However, this implies that at this step we are adding only one new quadric to $Q_{i'''}$, which is a contradiction since it implies that either $yw $ or $zw$ is in $Q_{i'''-1}$, which would imply that $w\in L_{i'''-1}$.
\end{proof}

\begin{example}
We exhibit some examples of the  partitions counted in \Cref{thm:expexp} for $s=3$ and $n=2,3$. In the caption, we describe the associated partition of $[n]$ (labelled with variables $x,y,z$) with the chosen coloured element. 

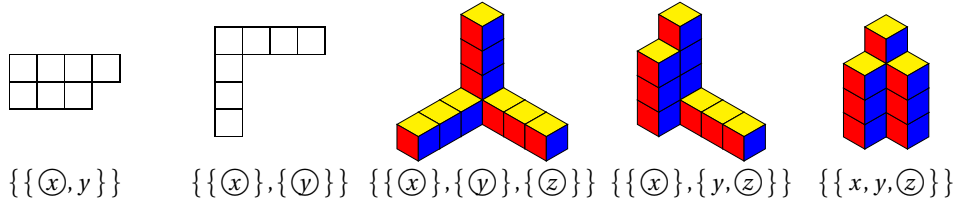
\begin{figure}[h!]
    \centering
  \scalebox{0.9}{\begin{tikzpicture}
    \node at (0,-1.5) {\small$\Set{\Set{\mbox{\textcircled{$x$}},y}}$}; 
    \node at (0,0) {\begin{ytableau}
\ &\   &\    &\    \\
\ &\   &\         \\   
\end{ytableau}};
    \node at (3,-1.5) {\small$\Set{\Set{\mbox{\textcircled{$x$}}},\Set{\mbox{\textcircled{$y$}}}}$}; 
\node at (3,0) {\begin{ytableau}
\ &\   &\    &\     \\
\  \\ \  \\  \   
\end{ytableau}};
    \node at (6.1,-1.5) {\small$\Set{\Set{\mbox{\textcircled{$x$}}},\Set{\mbox{\textcircled{$y$}}},\Set{\mbox{\textcircled{$z$}}}}$}; 
\node at (6.1,0) {\begin{tikzpicture}[scale=0.36] 
\planepartition{{4,1,1,1},{1},{1},{1}}  
\end{tikzpicture}};
    \node at (9.3,-1.5) {\small$\Set{\Set{\mbox{\textcircled{$x$}}},\Set{y,\mbox{\textcircled{$z$}}}}$}; 
\node at (9.3,0) {
\begin{tikzpicture}[scale=0.36] 
\planepartition{{4,1,1,1},{3}} 
\end{tikzpicture}};
    \node at (12,-1.5) {\small$\Set{\Set{x,y,\mbox{\textcircled{$z$}}}}$};
\node at (12,0) {\begin{tikzpicture}[scale=0.36] 
\planepartition{{4,3},{3}} 
\end{tikzpicture}}; 
\end{tikzpicture}}
    \caption{Examples of partitions from \Cref{thm:expexp} in the case $s=3$ and $n=2,3$.}
    \label{fig:esempio}
\end{figure}

Notice moreover that the monomial ideals
\begin{align*}
I_1&=\Set{x^3,y^3,x^2y,xy^2,zx,t}^{\perp} \subset \BC[x,y,z,t]\\
I_2&=\Set{x^3,y^3,x^2y,xy^2,zt}^{\perp} \subset \BC[x,y,z,t]   
\end{align*}
corresponding to partitions with Hilbert--Samuel function (1,4,4,4) are not of the form \eqref{eqn: forma delle partizioni step 3}, as they have  socle in degree 1, 2 and 3 (resp.~in degree 2 and 3).
\end{example}

\subsection{Hydral partitions and their counts}
In this subsection, we study the number of M-partitions with minimal number of boxes in degree 2 (cf.~\cref{prop:minq}), which 
comes in useful in view of the strategy we present in \Cref{subsec:algorithms}. We first explicitly derive closed formulas for the corresponding generating function when $n \in \set{2,3,4}$. In \Cref{subsetc:minimalquadgeneral}, we will present a general formula working in any dimension.

\begin{prop}
\label{prop:hydral-n=2}
Fix integers $s\geqslant 2$ and $m\geqslant 1$. Set $d=1+ 2(s-1)+m$ and let   $\overline \alpha_{s,m}$ be the number 
of  partitions $\lambda\in \mathrm{P}_d^2$ with $h_\lambda(2)=2$ and $\Soc(\lambda)\subset \lambda_{\geqslant s}$. Then, we have
\[
\overline \alpha_{s,m}=\overline \alpha_{3,m}=\alpha_{2,m}^2.
\]
Moreover, the generating function of these numbers is given by
\[
\sum_{m\geqslant 0}\overline{\alpha}_{s,m} t^m=\sum_{m\geqslant 1}\alpha_{2,m}^2t^m=\frac{t^2+t+2}{(1-t)(1-t^2)} t,
\]
showing that the left-hand side is a rational function.
\end{prop} 

\begin{proof}
The possible diagrams associated to linear partitions of the required form are classified according to $(1,1)\in \lambda$ or $(1,1)\notin\lambda$, see \Cref{fig:esempio2dim}  for the  two possible shapes (up to symmetry).

 \begin{figure}[h!]
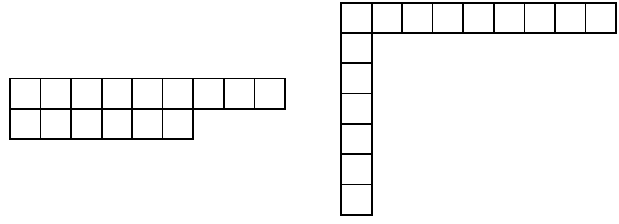

     \centering
     $\begin{matrix}
\begin{ytableau}
\ &\   &\    &\  &\  &\  &\ &\ &\   \\
\ &\ &\ &\ &\ &\  
\end{ytableau}
\end{matrix} \quad\quad \begin{matrix}
\begin{ytableau}
\ &\   &\    &\  &\  &\  &\ &\ &\   \\
\ \\\ \\\ \\\ \\\ \\\  
\end{ytableau}
\end{matrix} $
     \caption{On the left $(1,1) \in \lambda$, on the right $(1,1)\notin \lambda$.}
     \label{fig:esempio2dim}
 \end{figure}
 
Set
\begin{align*}
A_{1,s}^{(m)}&=\Set{\lambda\in \mathrm{P}_d^2  | h_\lambda(2)=2,  \Soc(\lambda)\in \lambda_{\geqslant s},\  (0,2)\notin \lambda},\\
A_{2,s}^{(m)}&=\Set{\lambda\in \mathrm{P}_d^2  | h_\lambda(2)=2,  \Soc(\lambda)\in \lambda_{\geqslant s},\  (1,1)\notin \lambda}.
\end{align*}
Then, we have
\[
\Lambda_{d,s}^{(2)}=2\lvert A_{1,s}^{(m)}\rvert + \lvert A_{2,s}^{(m)}\rvert.
\]
A simple computation yields
\[
\lvert A_{1,s}^{(m)}\rvert = \left\lfloor\frac{d}{2}\right\rfloor-s+1 ,\quad\quad \lvert A_{2,s}^{(m)}\rvert = d-2s.
\]
The final claim follows by  expanding the rational function and comparing the coefficients.
\end{proof}

\begin{remark}\label{rem:specialdim2} Notice that the number $\lvert A_{1,3}^{(m)}\rvert$ coincides with the number of partitions of $m-1$ into two parts (possibly empty), the first being greater or equal than the second. Indeed, if a linear partition $\lambda$ does not contain the point $(0,2)$, its socle is of the form
 \[
 \Soc(\lambda)=\Set{ (h,1),(h+j,0) },
 \]
 with $m-1=(h+j-2)+(h-2)$.

These combinatorial objects are known as \textit{headstrong combinations} of $m-1$ into $2$ parts, see \Cref{def: head comb}.

On the other hand, the computation of $|A_{2,3}^{(m)}|$ is easier and can be reduced to a counting problem of 0-dimensional partitions.
\end{remark}

\begin{definition}
\label{def: hydral}
We say that a partition $\lambda\in \mathrm{P}_d^k$ is  \emph{hydral} if $h_{\lambda}(1)=h_{\lambda}(2)=k $ and  $\Soc(\lambda)\subset \lambda_{\geqslant 3}$. Equivalently, a hydral partition is an M-partition of type $(k,k,m)$ for some $k,m\geqslant 1$. 
\end{definition}

\begin{prop}
\label{prop:dim3any}
Fix an integer $m\geqslant 1$ and set $d=1+6+m$.
The number  $\alpha_{3,m}^3$ of hydral partitions $\lambda\in \mathrm{P}_d^3$ is given by 
\[
\alpha_{3,m}^{3}=\begin{cases}
     1& \mbox{if }m=1\\ 
     6\left(m-1+\sum_{k=0}^{m-2}\left\lfloor\frac{k}{2}\right\rfloor \right)+3\left(1 + \left\lfloor\frac{m^2-4}{6}\right\rfloor\right) +\binom{m-1}{2} &\mbox{if }m\geqslant 2.
 \end{cases} 
 \]
Moreover, the generating function of these numbers is given by
\[
\sum_{m\geqslant 0}\alpha_{3,m}^{3} t^m=\frac{1+8t+6 t^2 +11  t^3 +3 t^4+2 t^5-t^6
}{(1-t)(1-t^2)(1-t^3)}t,
\]
showing that the left-hand side is a rational function.
\end{prop}

\begin{proof}
We consider the contribution to $\alpha_{3,m}^{3}$ from the 5 possible configurations of partitions according to their degree 2 boxes (up to symmetry). For each case, we display the corresponding $\lambda_{\leqslant 2}$.

   \underline{1.} $\lambda_{=2}=\Set{(2,0,0),(0,2,0),(0,0,2)}$.
        
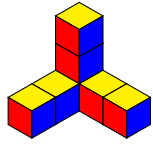
\begin{center}
            \begin{tikzpicture}[scale=0.36] 
\planepartition{{3,1,1},{1},{1}} 
\end{tikzpicture}
        \end{center}
        If $m=1,2$ this configuration does not contribute, while if $m\geqslant 3$ it contributes as 
        \[
        \binom{m-1}{2}.
        \]
\underline{2.} $\lambda_{=2}=\Set{(2,0,0),(0,2,0),(1,0,1)}$.
        
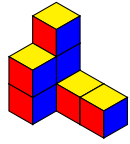
\begin{center}
            \begin{tikzpicture}[scale=0.36] 
\planepartition{{3,1,1},{2}} 
\end{tikzpicture}
\end{center}
        If $m=1$ this configuration does not contribute, while if $m\geqslant 2$ it contributes as  
        \[
        6\left(m-1+\sum_{k=0}^{m-2}\left\lfloor\frac{k}{2}\right\rfloor \right).
        \]
        \underline{3.} $\lambda_{=2}=\Set{(2,0,0),(1,1,0),(1,0,1)}$.
        
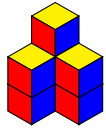
\begin{center}
            \begin{tikzpicture}[scale=0.36] 
\planepartition{{3,2},{2}} 
\end{tikzpicture}
        \end{center}
         If $m=1$ this configuration does not contribute, while if $m\geqslant 2$ it contributes as   
        \[
        3\left(1 + \left\lfloor\frac{m^2-4}{6}\right\rfloor\right).
        \]
    \underline{4.} 
        $\lambda_{=2}=\Set{(2,0,0),(1,1,0),(0,1,1)}$.
        \begin{center}
            \begin{tikzpicture}[scale=0.36] 
\planepartition{{3,2},{1,1}} 
\end{tikzpicture}
        \end{center}
       This configuration never contributes, as it fails the minimality of degree 2 boxes.\\
          \underline{5.} 
  $\lambda_{=2}=\Set{(1,0,1),(1,1,0),(0,1,1)}$.
        \begin{center}
            \begin{tikzpicture}[scale=0.36] 
\planepartition{{2,2},{2,1}} 
\end{tikzpicture}
        \end{center}
        If $m=1$ it contributes with 1, otherwise it does not contribute.

        The final claim follows by  expanding the rational function and checking its coefficients.
\end{proof}

\begin{remark}
As already observed  in \Cref{rem:specialdim2}, some configurations (the cases \underline{1.} and \underline{2.}) that contribute to $\alpha_{3,m}^3$ can be derived from  $\alpha_{n,m}^n$ for $n=1,2$.  
On the other hand, the cases \underline{3.} and \underline{5.} are purely three-dimensional. In particular, notice that \underline{3.} contributes as the number of partitions of $m-2$ in three ordered (possibly empty) parts, with  the first one being greater or equal than the others.

We will see in \Cref{subsetc:minimalquadgeneral} that these cases cover essentially all the possibilities in any dimension.
\end{remark}

\begin{example}
\label{ex:dim4}
Fix an integer $m\geqslant 1$ and set $d=1+8+m$. The number  $\alpha_{4,m}^4$ of hydral partitions $\lambda\in \mathrm{P}_d^4$ has the following generating  series 
\[
\sum_{m\geqslant 0}\alpha_{4,m}^4t^m=\frac{16+40t+57t^2+79t^3+81t^4+46t^5+37t^6+7t^7+t^8-4t^9}{(1-t^2)(1-t^3)(1-t^4)}t^2.
\] 
This formula is a direct application  of \Cref{cor: count hydral}.
\end{example}

\subsection{Headstrong M-partitions and their count}
\label{subsetc:minimalquadgeneral} 
We  prove  a closed formula for the number $\alpha_{n,m}^n$ of  hydral partitions
for all $n,m\geqslant 1$ in terms of the combinatorics of \emph{headstrong compositions}.

\begin{definition}
\label{def: head comb}
Fix $m,n\in \BZ_{\geqslant 1}$. A \emph{headstrong composition} of $m$ into $n$ parts is a combination (ordered partition) of $m$ into $n$ parts,\footnote{For our convenience, we choose our definition consistently with \cite[\href{https://oeis.org/A156041}{A156041}]{oeis}. Beware that in the literature, e.g.~\cite{Galgross}, some authors define headstrong compositions allowing only strictly positive $\alpha_i$. Nevertheless, the two associated generating functions are related to each other in an elementary way, cf.~\Cref{lemma:anydimcount}.} with the first part greater or equal than all the other parts. In other words,  it is an $n$-tuple $(\alpha_i)_{i=1}^n $ of nonnegative integers such that
    \begin{itemize}
        \item $\sum_{i=1}^n \alpha_i=m$,
        \item $\alpha_1\geqslant \alpha_i$ for all $i=1, \dots, n$.
    \end{itemize}
\end{definition}

\begin{definition}
Let $n\geqslant 2$ and $m\geqslant n-1$ be two positive integers.
An M-partitions $\lambda$ of type $(n,n,m)$ is called \emph{headstrong} if there exists a headstrong combination 
   $(\alpha_i)_{i=1}^n $ of $m-n+1$ into $n$ parts and a permutation $\sigma\colon [n]\to [n]$ such that  
\[
\lambda=\left(\Set{(\alpha_1+2)\cdot \mathbf{e}_{\sigma(1)}}\amalg \Set{\mathbf{e}_{\sigma(i)}+(\alpha_i+2)\cdot \mathbf{e}_{\sigma(1)}|i=2,\ldots,n} \right)^\perp\subset\BN^n,
\] 
where $\mathbf{e}_i$ denotes the $i$-th canonical generator of $\BN^n$.
\end{definition}
 
\begin{example}
An example of headstrong combinations of $10 $ into $4$ parts are
    \[
    \alpha=(4,1,2,3 )\mbox{ or }\beta=(5,1,3,1).
    \]
The permutation $\sigma=(12)(34) \in \mathfrak S_4$ induces two headstrong partitions $\pi_1,\pi_2$, which we display in monomial notation
    \[
\pi_1=    \Set{x_2^{6},x_1x_2^{3},x_4x_2^{4},x_3x_2^{5} }^\perp,\quad
\pi_2=    \Set{x_2^{7},x_1x_2^{3},x_4x_2^{5},x_3 x_2^{3}}^\perp.
    \] 
     \Cref{fig:headstrong}  depicts two example of headstrong M-partitions partitions of type $(3,3,4)$ and $(4,4,16)$ and compatible with the headstrong combinations $(3,0,1)$ and $(6,2,1,4)$. Two possible realisations are
\[
\Set{x_1^5,x_2x_1^2,x_3x_1^3}^\perp \mbox{ and } \Set{x_1^8,x_2x_1^3,x_3x_1^4,x_4x_1^6}^\perp.
\]
Other examples of headstrong partitions can be found in the first and last pictures of \Cref{fig:esempio} and in the left picture of \Cref{fig:esempio2dim}. 
    \begin{figure}[ht]
        \centering
        
            \begin{tikzpicture}[scale=0.36] 
\planepartition{{6,4},{3}}

\end{tikzpicture} 
 \hspace{2cm} 
   \begin{tikzpicture}[scale=0.7]
  \pgftransformcm{-0.707}{-0.4}{0.707}{-0.4}{\pgfpoint{0cm}{0cm}} 

  \coordinate (A1) at (0,0,0);
  \coordinate (B1) at (2,0,0);
  \coordinate (C1) at (2,2,0);
  \coordinate (D1) at (0,2,0);
  \coordinate (E1) at (0,0,4);
  \coordinate (F1) at (2,0,4);
  \coordinate (G1) at (2,2,4);
  \coordinate (H1) at (0,2,4);

  \coordinate (A2) at (0,-2,0);
  \coordinate (B2) at (2,-2,0);
  \coordinate (C2) at (2,0,0);
  \coordinate (D2) at (0,0,0);
  \coordinate (E2) at (0,-2,4);
  \coordinate (F2) at (2,-2,4);
  \coordinate (G2) at (2,0,4);
  \coordinate (H2) at (0,0,4);

  \draw[thick, fill=red!10, opacity=0.5] (A2) -- (B2) -- (C2) -- (D2) -- cycle;
  \draw[thick, fill=red!10, opacity=0.5] (A2) -- (B2) -- (F2) -- (E2) -- cycle;
  \draw[thick, fill=red!10, opacity=0.5] (A2) -- (E2) -- (H2) -- (D2) -- cycle;
  \draw[thick, fill=red!10, opacity=0.5] (E2) -- (F2) -- (G2) -- (H2) -- cycle;
  \draw[thick, fill=red!10, opacity=0.5] (B2) -- (F2) -- (G2) -- (C2) -- cycle;
  \draw[thick, fill=red!10, opacity=0.5] (D2) -- (C2) -- (G2) -- (H2) -- cycle;
  \node[align=center, font=\bfseries\Huge ] at (1.5,-1.7,1) {7}; 

  \coordinate (A3) at (-2,0,0);
  \coordinate (B3) at (0,0,0);
  \coordinate (C3) at (0,2,0);
  \coordinate (D3) at (-2,2,0);
  \coordinate (E3) at (-2,0,4);
  \coordinate (F3) at (0,0,4);
  \coordinate (G3) at (0,2,4);
  \coordinate (H3) at (-2,2,4);

  \draw[thick, fill=green!10, opacity=0.5] (A3) -- (B3) -- (F3) -- (E3) -- cycle;
  \draw[thick, fill=green!10, opacity=0.5] (A3) -- (B3) -- (C3) -- (D3) -- cycle;
  \draw[thick, fill=green!10, opacity=0.5] (A3) -- (E3) -- (H3) -- (D3) -- cycle;
  \draw[thick, fill=green!10, opacity=0.5] (E3) -- (F3) -- (G3) -- (H3) -- cycle;
  \draw[thick, fill=green!10, opacity=0.5] (B3) -- (F3) -- (G3) -- (C3) -- cycle;
  \draw[thick, fill=green!10, opacity=0.5] (D3) -- (C3) -- (G3) -- (H3) -- cycle;
  \node[align=center, font=\bfseries\Huge ] at (-1.7,1.5,1) {4}; 

  \coordinate (A4) at (0,0,4);
  \coordinate (B4) at (2,0,4);
  \coordinate (C4) at (2,2,4);
  \coordinate (D4) at (0,2,4);
  \coordinate (E4) at (0,0,8);
  \coordinate (F4) at (2,0,8);
  \coordinate (G4) at (2,2,8);
  \coordinate (H4) at (0,2,8);

  \draw[thick, fill=blue!10, opacity=0.5] (A1) -- (B1) -- (F1) -- (E1) -- cycle; 
  \draw[thick, fill=blue!10, opacity=0.5] (A1) -- (B1) -- (C1) -- (D1) -- cycle; 
  \draw[thick, fill=blue!10, opacity=0.5] (A1) -- (E1) -- (H1) -- (D1) -- cycle; 
  \draw[thick, fill=blue!10, opacity=0.5] (E1) -- (F1) -- (G1) -- (H1) -- cycle; 
  \draw[thick, fill=blue!10, opacity=0.5] (B1) -- (F1) -- (G1) -- (C1) -- cycle; 
  \draw[thick, fill=blue!10, opacity=0.5] (D1) -- (C1) -- (G1) -- (H1) -- cycle; 
  \node[align=center, font=\bfseries\Huge ] at (1.5,1.5,1) {9};
  \draw[thick, fill=yellow!10, opacity=0.5] (A4) -- (B4) -- (F4) -- (E4) -- cycle;
  \draw[thick, fill=yellow!10, opacity=0.5] (A4) -- (B4) -- (C4) -- (D4) -- cycle;
  \draw[thick, fill=yellow!10, opacity=0.5] (A4) -- (E4) -- (H4) -- (D4) -- cycle;
  \draw[thick, fill=yellow!10, opacity=0.5] (E4) -- (F4) -- (G4) -- (H4) -- cycle;
  \draw[thick, fill=yellow!10, opacity=0.5] (B4) -- (F4) -- (G4) -- (C4) -- cycle;
  \draw[thick, fill=yellow!10, opacity=0.5] (D4) -- (C4) -- (G4) -- (H4) -- cycle;
  \node[align=center, font=\bfseries\Huge ] at (1.8,1.8,6) {5};
\end{tikzpicture}
        \caption{Examples of planar and (a projection of) a solid headstrongs partitions.}
        \label{fig:headstrong}
    \end{figure}
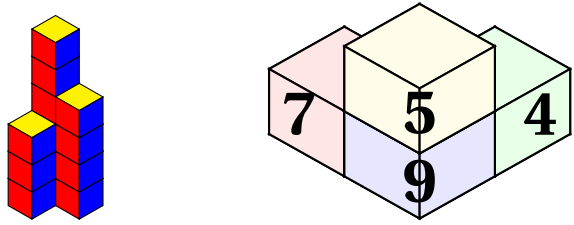
\end{example} 

\begin{remark}\label{rem: perm headstrong}
    If $\lambda $ is a headstrong M-partition of type $(n,n,m)$, then there exists a permutation $\sigma \in \mathfrak S_n$ such that
\[ 
\lambda_{=2}=\Set{\mathbf{e}_{\sigma(1)}+\mathbf{e}_{{\sigma(i)}}|i=1,\ldots,n}\subset\BN^n. 
\]
In \Cref{lemma:anydimcount} we will compute the number of headstrong M-partitions whose `head' corresponds to the first variable.
\end{remark}

\begin{lemma}
\label{lemma:anydimcount}  
Fix two integers $n\geqslant 1$ and $m\geqslant 0$. Denote by $\delta_{n,m}$ the number of M-partitions $\lambda$ of type $(n,n,m)$ such that 
    \begin{align}\label{eqn: condition star}
\lambda_{=2}=\Set{\mathbf{e}_{1}+\mathbf{e}_{{i}}|i=1,\ldots,n}\subset\BN^n.
\end{align}
Then, the generating function
\[
\Phi_n(t)=\sum_{m\geqslant 0}\delta_{n,m}t^{m}
\]
satisfies the identity
\begin{align*}
\Phi_n(t)=
\begin{cases}
\displaystyle\sum_{i\geqslant 1}t^i&\mbox{ for }n=1,\\ \\
\displaystyle\sum_{i\geqslant 1}t^{i+n-2}\left(\frac{1-t^i}{1-t}\right)^{n-1}&\mbox{ for }n\geqslant 2.
\end{cases}
\end{align*} 
In particular, the number of headstrong M-partitions of type $(n,n,m)$ is $n\cdot \delta_{n,m}$.
\end{lemma}

\begin{proof} The case $n=1$ is a direct check.  Let $n\geqslant 2$. Denote by $\Delta_{n,m}$ the collection of headstrong M-partitions  of type $(n,n,m)$ satisfying \eqref{eqn: condition star}. 
    By condition \eqref{eqn: condition star}, if $\lambda\in \Delta_{n,m}$ then
    \[  \Set{\mathbf{e}_{1}+\mathbf{e}_{{i}}|i=2,\ldots,n}\subset \lambda_{=3}.
    \]
    As a consequence, we have $\Delta_{n,m}=\varnothing$, for $m<n-1$. Let us now assume $m\geqslant n-1$ and let us denote by $\mathsf{HS} _{n,m-n+1}$ the set of headstrong compositions of $m-n+1$ into $n$ parts. There is a bijective correspondence
\[
\begin{tikzcd}[row sep =tiny]
    \mathsf{HS} _{n,m-n+1 }\arrow[r]&\Delta_{n,m}\\
     (\alpha_1,\ldots,\alpha_n) \arrow[r,mapsto] & \left(\Set{(\alpha_1+2)\cdot \mathbf{e}_{1}}\amalg \Set{\mathbf{e}_{i}+(\alpha_i+2)\cdot \mathbf{e}_{1}|i=2,\ldots,n} \right)^\perp.
\end{tikzcd}
\] 
By the proof of \cite[Thm.~20]{Galgross}, the generating function of the number of headstrong combinations of $m$ into $n$ parts is given by
\begin{align*}
    \sum_{m\geqslant 1}\lvert  \mathsf{HS} _{n,m }\rvert t^m=\sum_{i\geqslant 1}t^{i-1}\left(\frac{1-t^i}{1-t}\right)^{n-1}.
\end{align*}
Combining this identity with the bijective correspondence above yields the result.

The final claim follows from the fact that  there are $n$ choices for the head of the partition.
\end{proof}

\begin{remark}
\label{rem:rationhydral}
The generating function in \Cref{lemma:anydimcount} is rational: a straightforward manipulation yields the equality 
\[ 
\Phi_n(t) = \frac{t^{n-2}}{(1-t)^{n-1}}\sum_{k=1}^n\binom{n-1}{k-1}(-1)^{k+1}\frac{t^k}{1-t^k},
\]
see also \cite[Sec.~I.3]{FS_book} for some related rational generating functions.
\end{remark}
 
In the next result, we give a complete characterisation of hydral partitions, in terms of headstrong M-partitions.
\begin{theorem}\label{thm: hydral structure}
   Let $\lambda$ be a hydral partition of type $(k,k,m)$, for $k,m\geqslant 1$. Then there exist
   \begin{itemize}
       \item a partition of  the set $[k]=\Set{1, \dots, k}$
       \begin{equation}\label{eq:partitionhydral}
            [k]=\coprod_{i=1}^r A_i\amalg \coprod_{i=1}^s B_i \amalg \coprod_{i=1}^t C_i,
       \end{equation}
       such that 
       \begin{itemize}
           \item $ A_i =\Set{a_i}$ is a singleton for all $i=1, \dots, r$,
        \item    $|B_i|\geqslant 2$ for all $i=1, \dots, s$,
           \item $|C_i|=3$ for all $i=1, \dots, t$,
           \item $\widetilde{m}:=m- \left(r+t-s+\sum_{i=1}^s |B_i|\right)\geqslant 0,$
       \end{itemize}
       \item two tuples $(m_i)^r_{i=1}, (n_i)^s_{i=1}$ of nonnegative integers such that 
       \[
       \widetilde{m}=\sum_{i=1}^rm_i +\sum_{i=1}^s n_i, 
       \]
       \item for each $i=1, \dots, s$, a headstrong M-partition $\pi_i$ of type $(|B_i|, |B_i|, n_i+|B_i|-1)$,
   \end{itemize}
such that, in monomial notation, we have
\[
\lambda=\bigcup_{i=1}^r \Set{x_{a_i}^{m_i+3}}^\perp \cup\bigcup_{i=1}^s \pi_i \cup \bigcup_{i=1}^t\Set{x_ax_bx_c|C_i=\Set{a,b,c}}^\perp.
\]
\end{theorem}

\begin{proof}
Let $\lambda$ be a hydral partition of type $(k,k,m)$, and denote by $\overline{\lambda}=\lambda_{\leqslant 3}$. We divide the proof into several steps. To ease the notation, we denote the elements of $\lambda, \overline{\lambda}$ in terms of monomials, cf.~\Cref{sec: monomial ideals}.

\underline{Step I.}
Let $xyz\in \lambda_{=3} $ be a square-free cubic monomial and $\mathfrak{n} \in \lambda_{=3}$ any other cubic monomial. We claim that none of the three different variables $x,y,z$ appear in $\mathfrak{n}$. Assume by contradiction that there exists at least one variable in common. We prove the case when $\mathfrak{n}=xyz'$ is square-free with $z\neq z'$; the other cases will follow by an analogous analysis.

We proceed as in the proof of \Cref{thm:expexp}, `building' the partition $\overline{\lambda}$ by adding one by one the socle elements, and looking at their contributions to the Hilbert--Samuel function of $\overline{\lambda}$. Let 
  \begin{align*}
      W_0&=\varnothing,\\
      W_1&=W_0\cup \set{xyz}.
  \end{align*}
The  cubic $xyz$ contributes to the Hilbert--Samuel function with a vector $(1,3,3,1)$, i.e.~with one degree zero element, three linear terms, three quadrics and one cubic. Set now 
\[
W_2=W_1\cup \set{xyz'}.
\] 
The cubic $xyz'$ contributes  to the Hilbert--Samuel function with a vector $(0,1,2,1)$, for a total contribution of $(1,4,5,2)$. Since the Hilbert--Samuel function of $\overline{\lambda}$ is $(1,k,k,k')$ for some $k'\leqslant m$, there must exists a `step' $W_i$ where we add a cubic monomial $\mathfrak{m}$ which contributes with more variables than quadratic monomials. The only possibilities are that $\mathfrak{m}$ contributes with $(0,2,1,1)$ or $(0,3,2,1)$. If the contributing vector is $(0,2,1,1)$, it means that $\mathfrak{m}$ contributes with exactly two new variables $a,b$. Therefore, up to relabelling the variables,  its possibilities are
  \[
  \mathfrak{m}\in\set{a^2b, \, abc},
  \]
where $c$ is a variable already added at a previous step. In both cases, $\mathfrak{m}$ contributes with 2 new quadrics, which is a contradiction. The case of contributing vector $(0,3,2,1) $ follows from an analogous reasoning. 

If we denote by 
\[
\set{x_iy_iz_i| i=1, \dots, t}
\]
the set of square-free cubics in $\lambda$, we set $C_i=\set{x_i, y_i, z_i}$.

 \underline{Step II.}
By Step I,  if $\lambda$ contains a square-free monomial $xyz$, then the only monomials it contains with variables $x,y,z$ are precisely all the partial derivatives of $xyz$. Therefore, without loss of generality we can suppose that $\lambda $ contains no square-free cubic monomials. Denote  the Hilbert--Samuel function of $\overline{\lambda}$ by $(1,k,k,k')$. 

We claim that $k'\leqslant k$. Assume by contradiction that $k'>k$. Then it means that there exists in $\lambda$ at least one cubic of the form $x^2y$ with $x\neq y$. Let $x_1^3, \dots, x_\ell^3$ be all the pure cubic monomial in $\lambda$. As in Step I, set
 \begin{align*}
     W_0&=\varnothing,\\
     W_i&=W_{i-1}\cup \set{x_i^3}, \quad i=1, \dots, \ell.
 \end{align*}
 The first step contributes to the total Hilbert--Samuel function of $\overline{\lambda}$ with $(1,1,1,1)$ while all the others with $(0,1,1,1)$, for a total of $(1,\ell, \ell, \ell)$. Set 
\begin{align*}
    W_{\ell+1}=W_\ell\cup \set{x^2y}.
\end{align*} 
The possible vector contribution of this step are
\begin{align*}
    &(0,2,2,1)\\
    &(0,1,2,1)\\
    &(0,1,1,1)\\
    &(0,0,1,1).
\end{align*}
A  case by case analysis as in Step I shows that none of these situations  occurs.

  \underline{Step III.}
Consider cubic monomials $\mathfrak{m}_1, \dots, \mathfrak{m}_{k'-k}$ such that
\[\widetilde{\lambda}= \overline{\lambda}\cup\set{\mathfrak{m}_1, \dots, \mathfrak{m}_{k'-k}},
\]
is a partition with Hilbert--Samuel function $(1,k,k,k)$.  In particular,  $\mathfrak{m}_i$  belongs to the socle of  $\widetilde{\lambda} $, for all $i=1,\ldots,k'-k$. By the proof of \Cref{thm:expexp}, there exists a partition of the set of variables 
\begin{align*}
    [k]=\coprod_{i=1}^r A_i\amalg \coprod_{i=1}^s B_i,
\end{align*}
where $A_i=\set{a_i}$ and $|B_i|\geqslant 2$, with distinguished elements $b_i\in B_i$, such that
\begin{align*}
    \widetilde{\lambda}=\bigcup_{i=1}^r \Set{x_{a_i}^{3}}^\perp \cup\bigcup_{i=1}^s  \Set{x_{b_i}^{2}x_b  |  b\in B_i}^\perp.
\end{align*}
In particular, this implies that the monomials we possibly added can only be of the form $ x_{b_i}^3 $.

\underline{Step IV.} For a general hydral partition $\lambda$ of type $(k,k,m)$, we define the sets $C_\bullet$ as in Step I, and the sets $A_\bullet, B_\bullet$ as in Step III among the remaining variables. The sets $A_\bullet, B_\bullet, C_\bullet$ completely describe the partition $\lambda_{\leqslant 2}$. Denote by $\lambda_{A_\bullet}$ (resp. $\lambda_{B_\bullet}, \lambda_{C_\bullet}$) the projection of the partition $\lambda$ to the variables in $A_i$ (resp. in $B_i, C_i$). Clearly, $\lambda$ is reconstructed by its projections $\lambda_{A_\bullet}, \lambda_{B_\bullet}, \lambda_{C_\bullet}$.

By the  description of the elements in $\lambda$ of degree 3 in Step III, we  have that
\begin{align*}
    \lambda_{A_i}&=\set{x_{a_i}^{m_i+3}}^{\perp}, \quad \mbox{for some } m_i\geqslant 0,\\
    \lambda_{C_i}&= \set{x_iy_iz_i}^\perp,
\end{align*}
and that $\lambda_{B_i}$ is a headstrong M-partition of type $(|B_i|, |B_i|, n_i+|B_i|-1)$ for some $n_i$. The fact that $\widetilde{m}\geqslant 0$ follows by a simple counting of the number of  boxes in $\lambda$.
\end{proof}

\subsection{Proof of \texorpdfstring{\Cref{thm: intro hydral}}{}}
In this subsection we prove \Cref{thm: intro hydral} from the introduction, see \Cref{cor: count hydral}.
To this end, we need to introduce some auxiliary functions that will appear in the closed formula counting the number of hydral partitions. 

\begin{notation}
\label{not-for-thm-C}
 Recall from \Cref{notation:linear-partitions} that a linear partition $\lambda \in \mathrm P_d^2$ of a positive integer $d$ can be represented as  
\[
\lambda = (1^{\alpha_1}2^{\alpha_2}3^{\alpha_3}\cdots d^{\alpha_d}),
\]
where $\sum_{i=1}^d i\alpha_i=d$, and recall that $\lvert \aut(\lambda)\rvert =\prod_{i=1}^d\alpha_i!$. The integers $i\leqslant d$ such that $\alpha_i \neq 0$ form a decreasing tuple $\lambda_1'>\cdots>\lambda'_{s'}$ of positive integers, each coming with its own multiplicity $\alpha_{\lambda'_k}$ for $k=1,\ldots,s'$. We set $s_k = \alpha_{\lambda'_k}$ for simplicity, so that, for instance, $\lvert \aut(\lambda)\rvert =\prod_{k=1}^{s'}s_k!$.
After setting $s=\sum_{k=1}^{s'} s_k$, we can also represent $\lambda \vdash d$ as
\[
\lambda=(\lambda_1, \dots, \lambda_s), \qquad \lambda_1\geqslant \lambda_2 \geqslant \cdots \geqslant \lambda_s > 0,
\]
where 
\begin{align*}
\lambda_1=\cdots=\lambda_{s_1}&=\lambda_1'\\
\lambda_{s_1+1}=\cdots =\lambda_{s_1+s_2}&=\lambda_2'\\
\lambda_{s-s_{s'}+1}=\dots=\lambda_{s}&=\lambda'_{s'}.
\end{align*}
Define the following auxiliary functions:
\begin{itemize}
\item [(i)] For arbitrary $\lambda \vdash d$, define
\[
\mathsf t(\lambda)=\left|\Set{i=1, \dots, s|\lambda_i=3}\right|.
\]
This is the multiplicity $\alpha_3$.
\item [(ii)] if $\lambda$ satisfies $\mathsf t(\lambda)>0$ (i.e.~$\alpha_3 > 0$), define 
\[
\mathsf s(\lambda) = (1^{\alpha_1}2^{\alpha_2}3^{\alpha_3-1}\cdots d^{\alpha_d})\vdash d-3.
\]
\item [(iii)] For arbitrary $\lambda \vdash d$, define
\[
\displaystyle \mathsf f(\lambda) =\frac{1}{\lvert\aut(\lambda)\rvert}\prod_{j=1}^s\binom{d-\sum_{i=1}^{j-1}\lambda_i}{\lambda_{j}}\lambda_j.
\]
\item [(iv)] For arbitrary $\lambda \vdash d$, define
\[
\mathsf u(\lambda)=(1^{\alpha_1} 2^{\alpha_2}\cdots \lambda_1^{\alpha_{\lambda_1}-1})\vdash d-\lambda_1.
\]
\item [(v)] For arbitrary $\lambda \vdash d$, define
\[
\mathsf r(\lambda)=\sum_{i=1}^s\max(\lambda_i-1,1).
\]
\item [(vi)] for $m\geqslant \mathsf r(\lambda)$, define  
\[
\displaystyle \mu(\lambda,m)=\sum_{i_1=0}^{m-\mathsf r(\mathsf u(\lambda))}\ \sum_{i_2=0}^{m-i_1-\mathsf r(\mathsf u^2(\lambda))}\ \sum_{i_3=0}^{m-i_1-i_2-\mathsf r(\mathsf u^3(\lambda))}\  \cdots\  \sum_{i_{s-1}=0}^{m-\sum_{j=1}^{s-2}i_j-\mathsf r(\mathsf u^{s-1}(\lambda))}\prod_{j=1}^s\delta_{\lambda_j,i_j},
\]
where $i_s= m-\sum_{j=1}^{s-1}i_j $, $\mathsf u^\ell=\mathsf u\circ \dots \circ \mathsf u$ ($\ell$ times iteration) and $\delta_{n,m}$ are defined in \Cref{lemma:anydimcount}. For our convenience, we set $\mu(\lambda,m)=0$ for $m < \mathsf r (\lambda)$.
\end{itemize} 
\end{notation}

\begin{lemma}
\label{lemma:techfinal}
Let $\lambda$ be a linear partition and define the series
\[
\Psi_\lambda(t)= \sum_{m\geqslant 0}\mu(\lambda,m)t^m.
\]
Then there is an identity
\[
\Psi_\lambda(t)=\prod_{i=1}^s \Phi_{\lambda_i}(t),
\]
where we adopt the conventions of \Cref{not-for-thm-C} for $\lambda$, and $\Phi_\bullet(t)$ was defined in \Cref{lemma:anydimcount}.
In particular, by \Cref{rem:rationhydral}, the series  $\Psi_\lambda(t)$ is rational.
\end{lemma}

\begin{proof}
    We prove the statement by induction on the length $\ell(\lambda)\geqslant 1$. The base case is trivial.  Fix now some linear partition $\lambda\vdash n$. Then we have
\begin{align*} 
\Psi_\lambda(t)=&\sum_{m\geqslant 0}\displaystyle \mu(\lambda,m)t^m\\
=&\sum_{m\geqslant 0}\sum_{i_1=0}^{m-\mathsf r(\mathsf u(\lambda))}\ \sum_{i_2=0}^{m-i_1-\mathsf r(\mathsf u^2(\lambda))}\  
\ \cdots\  \sum_{i_{s-1}=0}^{m-\sum_{j=1}^{s-2}i_j-\mathsf r(\mathsf u^{s-1}(\lambda))}\prod_{j=1}^s\delta_{\lambda_j,i_j}t^m\\
=&\sum_{m\geqslant 0}\sum_{i_1=\mathsf r(\mathsf u(\lambda))}^{m}\delta_{\lambda_1,i_1-\mathsf r(\mathsf u(\lambda))} \sum_{i_2=0}^{m+\mathsf r(\mathsf u(\lambda))-i_1-\mathsf r(\mathsf u^2(\lambda))}\  
\  \cdots\  \sum_{i_{s-1}=0}^{m+\mathsf r(\mathsf u(\lambda))-\sum_{j=1}^{s-2}i_j-\mathsf r(\mathsf u^{s-1}(\lambda))}\prod_{j=2}^s\delta_{\lambda_j,i_j}t^m\\
=&\sum_{m\geqslant \mathsf r(\mathsf u(\lambda))}\sum_{i_1=\mathsf r(\mathsf u(\lambda))}^{m}\delta_{\lambda_1,i_1-\mathsf r(\mathsf u(\lambda))} \sum_{i_2=0}^{m+\mathsf r(\mathsf u(\lambda))-i_1-\mathsf r(\mathsf u^2(\lambda))}\  
\  \cdots\  \sum_{i_{s-1}=0}^{m+\mathsf r(\mathsf u(\lambda))-\sum_{j=1}^{s-2}i_j-\mathsf r(\mathsf u^{s-1}(\lambda))}\prod_{j=2}^s\delta_{\lambda_j,i_j}t^m\\
=&\sum_{i_1\geqslant \mathsf r(\mathsf u(\lambda))}\sum_{m\geqslant i_1}\delta_{\lambda_1,i_1-\mathsf r(\mathsf u(\lambda))} \sum_{i_2=0}^{m+\mathsf r(\mathsf u(\lambda))-i_1-\mathsf r(\mathsf u^2(\lambda))}\ 
\  \cdots\  \sum_{i_{s-1}=0}^{m+\mathsf r(\mathsf u(\lambda))-\sum_{j=1}^{s-2}i_j-\mathsf r(\mathsf u^{s-1}(\lambda))}\prod_{j=2}^s\delta_{\lambda_j,i_j}t^m\\
=&\sum_{i_1\geqslant \mathsf r(\mathsf u(\lambda))}\sum_{m\geqslant 0}\delta_{\lambda_1,i_1-\mathsf r(\mathsf u(\lambda))} \sum_{i_2=0}^{m+\mathsf r(\mathsf u(\lambda))-\mathsf r(\mathsf u^2(\lambda))}\  
\  \cdots\  \sum_{i_{s-1}=0}^{m+\mathsf r(\mathsf u(\lambda))-\sum_{j=2}^{s-2}i_j-\mathsf r(\mathsf u^{s-1}(\lambda))}\prod_{j=2}^s\delta_{\lambda_j,i_j}t^{m+i_1}\\
=&\sum_{i_1\geqslant  0}\delta_{\lambda_1,i_1 } t^{i_1 }\sum_{m\geqslant 0}\sum_{i_2=0}^{m+\mathsf r(\mathsf u(\lambda))-\mathsf r(\mathsf u^2(\lambda))}\  
\  \cdots\  \sum_{i_{s-1}=0}^{m+\mathsf r(\mathsf u(\lambda))-\sum_{j=2}^{s-2}i_j-\mathsf r(\mathsf u^{s-1}(\lambda))}\prod_{j=2}^s\delta_{\lambda_j,i_j}t^{m+\mathsf r(\mathsf u(\lambda))}\\
=&\sum_{i_1\geqslant 0}\delta_{\lambda_1,i_1 } t^{i_1}\sum_{m\geqslant \mathsf r(\mathsf u(\lambda))}\mu(\mathsf u(\lambda),m)t^{m}\\
=& \Phi_{{\lambda}_1}(t)      \Psi_{\mathsf u(\lambda)}(t)\\
=& \prod_{i=1}^s \Phi_{\lambda_i}(t),
     \end{align*} 
as required.
\end{proof}

The factorisation of $\Psi_\lambda(t)$ into the rational functions $\Phi_{\lambda_i}(t)$ reflects the fact that, by \Cref{thm: hydral structure}, headstrong combinations are the building blocks for hydral partitions.

We can now complete the proof of \Cref{thm: intro hydral}.

\begin{corollary}
\label{cor: count hydral}
Let $n, m \in \BZ_{\geqslant 1}$. Then 
\[
\alpha_{n,m}^n=\sum_{\lambda\vdash n}\sum_{i=0}^{\mathsf t(\lambda)}\binom{n}{3i}\frac{(3i)!}{6^ii!} \mathsf f(\mathsf s^i(\lambda))\mu(\mathsf s^i(\lambda),m-i) .
\] 
Moreover, there is an identity
\begin{equation}\label{eq:serialpha}
\sum_{m\geqslant 0} \alpha_{n,m}^nt^m=\sum_{\lambda\vdash n}\sum_{i=0}^{\mathsf t(\lambda)}\binom{n}{3i}\frac{(3i)!}{6^ii!} \mathsf f(\mathsf s^i(\lambda)) \Psi_{\mathsf{s}^i(\lambda)}(t)t^i.
\end{equation}
In particular, the left-hand side of \eqref{eq:serialpha} is a rational function.
\end{corollary}

\begin{proof}
By \Cref{thm: hydral structure} the number $\alpha_{n,m}^n$ is given by the number of possibilities for the combinatorial data needed to describe a hydral partition.

First we claim that, given a linear partition $\lambda $ of size $n$, the number of ways to partition the set $[n]$ as in \eqref{eq:partitionhydral} with $t=0$ is
\[
\mathsf f (\lambda) \mu(\lambda,m).
\]
In fact, we have that
\[
\mathsf{f}(\lambda)=\Set{(D_i,d_i)|[n]=\coprod_{i=1}^s D_i,\ |D_i|=\lambda_i\mbox{ and } d_i\in D_i}.
\] 
On the other hand, $\mu(\lambda,m)$ counts  the number of partitions $\pi$ that can be written as a union
\[
\pi=\bigcup_{i=1}^s\pi_i
\]
of headstrong partitions $\pi_i$ such that
\[
({\pi_i})_{=2}=\Set{\mathbf{e}_{\sigma_i(1)}+\mathbf{e}_{{\sigma_i(j)}}|j=1,\ldots,\lambda_i}\subset\BN^n, 
\] 
for some bijection $\sigma_i\colon [\lambda_i]\to D_i$ sending $1 $ to $d_i$.

Secondly, we take into account the presence of the sets of variables $C_i$ as in \eqref{eq:partitionhydral}.  For $j=0,\ldots,\mathsf t(\lambda)$, we  have $\binom{n}{3i}$ ways to  choose $3i$ variables among the $n$ given to define the sets $C_1, \dots, C_j$ of cardinality 3. Once choosen the $3i$ preferred variables, there are $\frac{(3i)!}{6^ii!}$ ways to form the sets $C_i$. Therefore, there are in  total,$\binom{n}{3i} \frac{(3i)!}{6^ii!}$ ways to choose the sets $C_1, \dots, C_j$ as in  \eqref{eq:partitionhydral}. Once fixed the sets $C_i$, there remain $m-j$ boxes of $\lambda$ to be distributed in the variables corrisponding to indices in  $A_\bullet, B_\bullet$. By the first part of this proof,  this can be done in 
     \[
     \mathsf f (\mathsf s^i(\lambda)) \mu( \mathsf s^i(\lambda),m-i) 
     \]
     ways. This proves the first claimed equality.

   A simple manipulation shows that \eqref{eq:serialpha}  holds. Its rationality is  deduced by combining \cref{lemma:techfinal} and \Cref{rem:rationhydral}.
\end{proof}

\begin{remark}
\label{rem:optimality}
The rationality statement in \Cref{cor: count hydral} should be considered `optimal', in the sense that counting M-partitions with nonminimal number of quadrics, i.e.~with $h_\lambda(2)=h_\lambda(1)+1 $,  yields manifestly nonrational generating functions in any dimension.

As an example, we computed explictly a closed  formula in the case $h_\lambda(1)=3,h_\lambda(2)=4 $. Due to the complexity of the formula we do not display it, but rather discuss how to obtain it. Let $\lambda$ be an M-partition of type $(3,4,m)$. Then, up to symmetry, the partition $\lambda_{\leqslant 2}$ is of one of the forms presented in \Cref{fig:possshapenonhydral} below.
\begin{figure}[ht]
    \centering
    \begin{tikzpicture}[scale=0.36] 
\planepartition{{2,2,1},{2,1}}
\node at (0,-3) {\underline{1.}};
\end{tikzpicture}\hspace{1cm}\begin{tikzpicture}[scale=0.36] 
\planepartition{{3,1,1},{2,1}} 
\node at (0,-3) {\underline{2.}};
\end{tikzpicture}\hspace{1cm}\begin{tikzpicture}[scale=0.36] 
\planepartition{{2,2,1},{1,1},{1}}
\node at (0,-3) {\underline{3.}};
\end{tikzpicture}\hspace{1cm}\begin{tikzpicture}[scale=0.36] 
\planepartition{{3,1,1},{1,1},{1}}
\node at (0,-3) {\underline{4.}};
\end{tikzpicture}
\caption{Possible shapes of $\lambda\in \mathrm{P}_8^3$ with $h_\lambda(1) = 3$ and $\ell(\lambda)=2$.}
\label{fig:possshapenonhydral}
\end{figure}
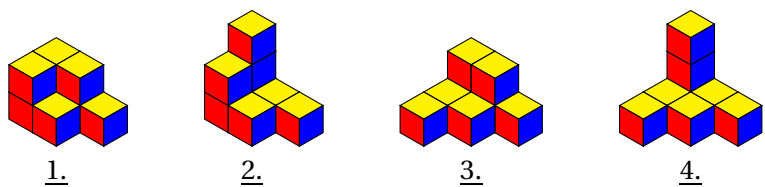
One can easily show that the configurations \underline{1.} and \underline{2.} contribute with rational generating functions that can be found in  \cite[\href{https://oeis.org/A008763}{A008763}, \href{https://oeis.org/A006918}{A006918}]{oeis}. On the other hand, configurations \underline{3.} and \underline{4.} do not.
In fact, as suggested by \Cref{fig:possshapenonhydral}, the counting of  M-partitions with shape as in  \underline{3.} and \underline{4.} is equivalent (up to a shift) to  the counting of linear partitions $\pi$ with $h_\pi(2)=3$, i.e. with the maximal number of quadrics. This counting clearly yields a nonrational generating series, as it differs  from the counting of all linear partitions (up to a finite shift) by the rational series in \Cref{prop:hydral-n=2}.
    \begin{table}[ht]
        \centering
        \begin{tabular}{c|ccccccccccccc }
        $m$     &1&2&3&4&5&6&7&8&9& 10&11&12&13   \\
        \hline
        $\alpha_{4,m}^3$ & 0& 18& 51& 126& 252& 474& 801& 1302& 2001& 3000& 4344& 6183& 8595   
        \end{tabular}
        \caption{Values of $\alpha_{4,m}^3$ for $m=1,\ldots,13$.}
        \label{tab:h1piu1}
    \end{table}
\end{remark}

\subsection{The computation of the numbers \texorpdfstring{$\alpha_{q,m}^k$}{}}\label{subsec:algorithms} 
 The naive approach to compute the numbers $\alpha_{\bullet,\bullet}^\bullet$ is to fix a \emph{degree} $d$, classify the possible admissible  subsets $S\subset \BN^k$ and their associated partitions $S^\perp$, and count only the correct ones. This naive approach turns out to be computationally too hard already for small values of $k$. 

We propose a new algorithic procedure to compute higher dimensional partitions exploiting apolarity of monomial ideals.
This approach allows us to control the socle type and hence to effectively construct M-partitions of given type $(k,q,m)$ and length $\ell$. 
 
Let us denote by $Q^k$ the set of quadrics
\[
Q^k=\Set{\boldit{x}\in \BN^k | \deg (\boldit{x}) =2}.
\]

We say that a subset  $U\subset Q^k$ is M-\emph{stable} if there exists an M-partition $\lambda$ such that $h_\lambda(1)=k$ and $\lambda_{=2}=U$. 

Let $q\geqslant 1$ be a positive integer and let $S_q^k$ be the set
\[ 
S_q^k=\Set{U\subset Q^k | U \mbox{ is M-stable, and } \lvert U\rvert=q}\big/\FS_k,
\] 
where $\FS_k$ is the symmetric group on $k$ letters.

\begin{example}
    For $k=2$ the set $\Set{(2,0)}$ is M-stable, while $\Set{(1,1)}$ is not. In fact, any M-partition containing $(1,1)$ also contains at least one among $(2,0)$ and $(0,2)$.
\end{example}

Fix an integer $\ell \geqslant 3$ and  $[U]\in S^k_q$. Define the partition
\[
\lambda_{U,\ell}=\bigcup_{\substack{ \ell(\lambda)=\ell,\\\lambda_{=2}=U} }\lambda.
\] 
Then, we have
\begin{equation}
\label{eq:alphabtt}
\alpha_{q,m,\ell}^k
=\sum_{[U]\in S_q^k}\lvert[U]\rvert\cdot \left\rvert
\Set{
\lambda\subset \lambda_{U,\ell} | 
\begin{array}{cc}
\lambda \textrm{ is an M-partition of} \\
\textrm{type }(k,q,m)
\textrm{ and }\ell(\lambda)=\ell
\end{array}
}
\right\rvert.
\end{equation}
Notice that, by definition of $\lambda_{U,\ell}$, an apolar partition $\lambda=S^\perp$ verifies $\lambda\subset\lambda_{U,\ell}$ if and only if $S\subset \lambda_{U,\ell}$. Moreover, in this case, we get $(S^\perp)_{=2}\subseteq U$. This observation, together with the  \eqref{eq:alphabtt}, provides an accessible set of possible socles for partitions contributing to the computation of $\alpha_{q,m,\ell}^k$. 

It is worth mentioning that another strong restriction is provided by Macaulay's bound for monomial ideals, see e.g.~\cite[Thm.~3.3]{INITALMONOMIAL}. Indeed, prescribing the number $q=h_\lambda(2)$ imposes constraints on the other values of the Hilbert--Samuel function $h_\lambda$, therefore  constraining the number of possible admissible subsets of $\lambda_{U,\ell}$.

\subsubsection{Algorithm}
The above strategy can be put into an algorithmic form to compute the numbers $\alpha_{q,m}^k$. To do so, one has to list explicitly all the possible socles of an M-partition of length $\ell$, i.e.~all admissible subsets of $\lambda_{U,\ell}$, for all $[U]\in S^k_q$, and for each of them to compute the Hilbert--Samuel function of the corresponding apolar partition. 
Via this procedure, one can  compute the numbers $\alpha_{q,m,\ell}^k$, which can be summed together to get the numbers $\alpha_{q,m}^k$ for every value of $k\leqslant q \leqslant \binom{k+1}{2}$, see  \Cref{prop:minq}. 

We have implemented this algorithm in Sage, Mathematica and Macaulay2  and computed enough values $c_e^k$ to recover the polynomials
$\mathsf c_x(t)$, for $x\leqslant 26$. Moreover, using the same strategy, we computed enough values $p_d^n$ to recover the polynomials $\mathsf h_d(t)$ for $d\leqslant 30$, thus extending the data provided in \cite{ThePartitionsProjectWebsite}. We have collected these polynomials in the ancillary file \ancillary.

\section{Conjectures}
\label{sec:conjectures}
\subsection{MacMahon's discrepancy}
\label{sec:MacMahon-discrepancy}

We conjectured in \cite[Sec.~9.4]{MOTIVEMOTIVEMOTIVE} a structural result on the discrepancy between the generating function of higher dimensional partitions and MacMahon's original guess. We present here a stronger version of \cite[Conj.~9.8]{MOTIVEMOTIVEMOTIVE}.

Define numbers $ \pi_d^{n},$ by
\[
 \sum_{d\geqslant 0}\pi_d^{n}t^d =\prod_{m>0} \left(1-t^m\right)^{-\binom{m+n-3}{n-2}}, 
\]
In \cite{Amanov-Yeliussizov} it has been proven that $\pi_d^n$ counts the number of $(n-1)$-dimensional partitions (of varying size) having corner-hook volume equal to $d$. 
According to \cite[p.~1098]{Atkin}, there are integers $\overline{y}^\bullet_\bullet$ such that
\begin{align}\label{eqn: ybarrrre}
    \pi_d^n  = \sum_{k=0}^{d-1}\binom{n}{k} \overline{y}_{d}^k.
\end{align}
Analogously to \Cref{prop: inversion}\ref{inversion-1}, the numbers $\overline{y}^\bullet_\bullet$ can be interpreted as a refinement of the numbers $\pi_\bullet^\bullet$.

\begin{remark}
We do not report here the definition of \emph{corner-hook volume} as it is beyond the scope of this paper, see \cite{Amanov-Yeliussizov}. Nevertheless, we give some example. In dimension $n\le 2$ corner-hook volume and size agree. Even though for $n=3$ we have $p_d^3=\pi_d^3$, the two numbers count different objects, see \Cref{fig:examplepim}.
    \begin{figure}[ht]
        \centering \begin{tabular}{c|c} 
        $\pi_3^3$ & $p_3^3$\\ 
        \begin{tikzpicture}[scale=0.37]
            \node at (0,0) {    \begin{tikzpicture}[scale=0.20] 
\planepartition{{3}}
        \end{tikzpicture}};
            \node at (3,0) {    \begin{tikzpicture}[scale=0.20] 
\planepartition{{1,1,1}}
        \end{tikzpicture}};
            \node at (6,0) {    \begin{tikzpicture}[scale=0.20] 
\planepartition{{1},{1},{1}}
        \end{tikzpicture}};
            \node at (9,0) {    \begin{tikzpicture}[scale=0.20] 
\planepartition{{2,1}}
        \end{tikzpicture}};
            \node at (12,0) {    \begin{tikzpicture}[scale=0.20] 
\planepartition{{2},{1}}
        \end{tikzpicture}};
            \node at (15,0) {    \begin{tikzpicture}[scale=0.20] 
\planepartition{{1,1},{1,1}}
        \end{tikzpicture}}; 
        \end{tikzpicture}&
              
        \begin{tikzpicture}[scale=0.37]
            \node at (0,0) {    \begin{tikzpicture}[scale=0.20] 
\planepartition{{3}}
        \end{tikzpicture}};
            \node at (3,0) {    \begin{tikzpicture}[scale=0.20] 
\planepartition{{1,1,1}}
        \end{tikzpicture}};
            \node at (6,0) {    \begin{tikzpicture}[scale=0.20] 
\planepartition{{1},{1},{1}}
        \end{tikzpicture}};
            \node at (9,0) {    \begin{tikzpicture}[scale=0.20] 
\planepartition{{2,1}}
        \end{tikzpicture}};
            \node at (12,0) {    \begin{tikzpicture}[scale=0.20] 
\planepartition{{2},{1}}
        \end{tikzpicture}};
            \node at (15,0) {    \begin{tikzpicture}[scale=0.20] 
\planepartition{{1,1},{1}}
        \end{tikzpicture}}; 
        \end{tikzpicture}\\
        $\pi_4^3$ & $p_4^3$\\ 
        \begin{tikzpicture}[scale=0.37]
            \node at (0,0) {    \begin{tikzpicture}[scale=0.20] 
\planepartition{{4}}
        \end{tikzpicture}};
            \node at (3,0) {    \begin{tikzpicture}[scale=0.20] 
\planepartition{{1,1,1,1}}
        \end{tikzpicture}};
            \node at (6,0) {    \begin{tikzpicture}[scale=0.20] 
\planepartition{{1},{1},{1},{1}}
        \end{tikzpicture}};
            \node at (9,0) {    \begin{tikzpicture}[scale=0.20] 
\planepartition{{3,1}}
        \end{tikzpicture}};
            \node at (12,0) {    \begin{tikzpicture}[scale=0.20] 
\planepartition{{3},{1}}
        \end{tikzpicture}}; 
        
            \node at (12,-3) {    \begin{tikzpicture}[scale=0.20] 
\planepartition{{1,1},{1}}
        \end{tikzpicture}};
            \node at (0,-3) {    \begin{tikzpicture}[scale=0.20] 
\planepartition{{2},{2}}
        \end{tikzpicture}};
            \node at (3,-3) {    \begin{tikzpicture}[scale=0.20] 
\planepartition{{2,2}}
        \end{tikzpicture}};
            \node at (6,-3) {    \begin{tikzpicture}[scale=0.20] 
\planepartition{{2},{1},{1}}
        \end{tikzpicture}};
            \node at (9,-3) {    \begin{tikzpicture}[scale=0.20] 
\planepartition{{2,1,1}}
        \end{tikzpicture}}; 
        
            \node at (3,-6) {    \begin{tikzpicture}[scale=0.20] 
\planepartition{{1,1,1},{1,1,1}}
        \end{tikzpicture}};
            \node at (6,-6) {    \begin{tikzpicture}[scale=0.20] 
\planepartition{{1,1},{1,1},{1,1}}
        \end{tikzpicture}};
            \node at (9,-6) {    \begin{tikzpicture}[scale=0.20] 
\planepartition{{2,1},{1,1}}
        \end{tikzpicture}}; 
        \end{tikzpicture} &    \begin{tikzpicture}[scale=0.37]
            \node at (0,0) {    \begin{tikzpicture}[scale=0.20] 
\planepartition{{4}}
        \end{tikzpicture}};
            \node at (3,0) {    \begin{tikzpicture}[scale=0.20] 
\planepartition{{1,1,1,1}}
        \end{tikzpicture}};
            \node at (6,0) {    \begin{tikzpicture}[scale=0.20] 
\planepartition{{1},{1},{1},{1}}
        \end{tikzpicture}};
            \node at (9,0) {    \begin{tikzpicture}[scale=0.20] 
\planepartition{{3,1}}
        \end{tikzpicture}};
            \node at (12,0) {    \begin{tikzpicture}[scale=0.20] 
\planepartition{{3},{1}}
        \end{tikzpicture}}; 
        
            \node at (12,-3) {    \begin{tikzpicture}[scale=0.20] 
\planepartition{{1,1},{1,1}}
        \end{tikzpicture}};
            \node at (0,-3) {    \begin{tikzpicture}[scale=0.20] 
\planepartition{{2},{2}}
        \end{tikzpicture}};
            \node at (3,-3) {    \begin{tikzpicture}[scale=0.20] 
\planepartition{{2,2}}
        \end{tikzpicture}};
            \node at (6,-3) {    \begin{tikzpicture}[scale=0.20] 
\planepartition{{2},{1},{1}}
        \end{tikzpicture}};
            \node at (9,-3) {    \begin{tikzpicture}[scale=0.20] 
\planepartition{{2,1,1}}
        \end{tikzpicture}};

            \node at (3,-6) {    \begin{tikzpicture}[scale=0.20] 
\planepartition{{1,1,1},{1}}
        \end{tikzpicture}};
            \node at (6,-6) {    \begin{tikzpicture}[scale=0.20] 
\planepartition{{1,1},{1},{1}}
        \end{tikzpicture}};
            \node at (9,-6) {    \begin{tikzpicture}[scale=0.20] 
\planepartition{{2,1},{1}}
        \end{tikzpicture}}; 
        \end{tikzpicture}\\
        $\pi_5^3$ & $p_5^3$\\ 
        \begin{tikzpicture}[scale=0.37]
            \node at (0.4,0) {    \begin{tikzpicture}[scale=0.20] 
\planepartition{{5}}
        \end{tikzpicture}};
            \node at (3,0) {    \begin{tikzpicture}[scale=0.20] 
\planepartition{{1,1,1,1,1}}
        \end{tikzpicture}};
            \node at (6,0) {    \begin{tikzpicture}[scale=0.20] 
\planepartition{{1},{1},{1},{1},{1}}
        \end{tikzpicture}};
            \node at (9,0) {    \begin{tikzpicture}[scale=0.20] 
\planepartition{{4,1}}
        \end{tikzpicture}};
            \node at (12,0) {    \begin{tikzpicture}[scale=0.20] 
\planepartition{{4},{1}}
        \end{tikzpicture}};
            \node at (15,0) {    \begin{tikzpicture}[scale=0.20] 
\planepartition{{3,1,1}}
        \end{tikzpicture}};
        
            \node at (0.4,-3) {    \begin{tikzpicture}[scale=0.20] 
\planepartition{{3},{1},{1}}
        \end{tikzpicture}};
            \node at (3,-3) {    \begin{tikzpicture}[scale=0.20] 
\planepartition{{3},{2}}
        \end{tikzpicture}};
            \node at (6,-3) {    \begin{tikzpicture}[scale=0.20] 
\planepartition{{3,2}}
        \end{tikzpicture}};
            \node at (9,-3) {    \begin{tikzpicture}[scale=0.20] 
\planepartition{{2},{1},{1},{1}}
        \end{tikzpicture}};
            \node at (12,-3) {    \begin{tikzpicture}[scale=0.20] 
\planepartition{{2,1,1,1}}
        \end{tikzpicture}};
            \node at (15,-3) {    \begin{tikzpicture}[scale=0.20] 
\planepartition{{2,2,1}}
        \end{tikzpicture}};
        
            \node at (0.4,-6) {    \begin{tikzpicture}[scale=0.20] 
\planepartition{{2},{2},{1}}
        \end{tikzpicture}};
            \node at (3,-6) {    \begin{tikzpicture}[scale=0.20] 
\planepartition{{1,1,1},{1}}
        \end{tikzpicture}};
            \node at (6,-6) {    \begin{tikzpicture}[scale=0.20] 
\planepartition{{1,1},{1},{1}}
        \end{tikzpicture}};
            \node at (9,-6) {    \begin{tikzpicture}[scale=0.20] 
\planepartition{{1,1,1,1},{1,1,1,1}}
        \end{tikzpicture}};
            \node at (12,-6) {    \begin{tikzpicture}[scale=0.20] 
\planepartition{{1,1},{1,1},{1,1},{1,1}}
        \end{tikzpicture}};
            \node at (15,-6) {    \begin{tikzpicture}[scale=0.20] 
\planepartition{{1,1,1},{1,1,1},{1,1,1}}
        \end{tikzpicture}};
        
            \node at (0.4,-9) {    \begin{tikzpicture}[scale=0.20] 
\planepartition{{3,1},{1,1}}
        \end{tikzpicture}};
            \node at (3,-9) {    \begin{tikzpicture}[scale=0.20] 
\planepartition{{2,1},{1,1},{1,1}}
        \end{tikzpicture}};
            \node at (6,-9) {    \begin{tikzpicture}[scale=0.20] 
\planepartition{{2,1,1},{1,1,1}}
        \end{tikzpicture}};
            \node at (9,-9) {    \begin{tikzpicture}[scale=0.20] 
\planepartition{{2,2},{1,1}}
        \end{tikzpicture}};
            \node at (12,-9) {    \begin{tikzpicture}[scale=0.20] 
\planepartition{{2,1},{2,1}}
        \end{tikzpicture}};
            \node at (15,-9) {    \begin{tikzpicture}[scale=0.20] 
\planepartition{{2,1},{1}}
        \end{tikzpicture}};
        \end{tikzpicture}&
              
        \begin{tikzpicture}[scale=0.37]
            \node at (0.4,0) {    \begin{tikzpicture}[scale=0.20] 
\planepartition{{5}}
        \end{tikzpicture}};
            \node at (3,0) {    \begin{tikzpicture}[scale=0.20] 
\planepartition{{1,1,1,1,1}}
        \end{tikzpicture}};
            \node at (6,0) {    \begin{tikzpicture}[scale=0.20] 
\planepartition{{1},{1},{1},{1},{1}}
        \end{tikzpicture}};
            \node at (9,0) {    \begin{tikzpicture}[scale=0.20] 
\planepartition{{4,1}}
        \end{tikzpicture}};
            \node at (12,0) {    \begin{tikzpicture}[scale=0.20] 
\planepartition{{4},{1}}
        \end{tikzpicture}};
            \node at (15,0) {    \begin{tikzpicture}[scale=0.20] 
\planepartition{{3,1,1}}
        \end{tikzpicture}};
        
            \node at (0.4,-3) {    \begin{tikzpicture}[scale=0.20] 
\planepartition{{3},{1},{1}}
        \end{tikzpicture}};
            \node at (3,-3) {    \begin{tikzpicture}[scale=0.20] 
\planepartition{{3},{2}}
        \end{tikzpicture}};
            \node at (6,-3) {    \begin{tikzpicture}[scale=0.20] 
\planepartition{{3,2}}
        \end{tikzpicture}};
            \node at (9,-3) {    \begin{tikzpicture}[scale=0.20] 
\planepartition{{2},{1},{1},{1}}
        \end{tikzpicture}};
            \node at (12,-3) {    \begin{tikzpicture}[scale=0.20] 
\planepartition{{2,1,1,1}}
        \end{tikzpicture}};
            \node at (15,-3) {    \begin{tikzpicture}[scale=0.20] 
\planepartition{{2,2,1}}
        \end{tikzpicture}};
        
            \node at (0.4,-6) {    \begin{tikzpicture}[scale=0.20] 
\planepartition{{2},{2},{1}}
        \end{tikzpicture}};
            \node at (3,-6) {    \begin{tikzpicture}[scale=0.20] 
\planepartition{{1,1,1,1},{1}}
        \end{tikzpicture}};
            \node at (6,-6) {    \begin{tikzpicture}[scale=0.20] 
\planepartition{{1,1},{1},{1},{1}}
        \end{tikzpicture}};
            \node at (9,-6) {    \begin{tikzpicture}[scale=0.20] 
\planepartition{{1,1,1},{1,1}}
        \end{tikzpicture}};
            \node at (12,-6) {    \begin{tikzpicture}[scale=0.20] 
\planepartition{{1,1},{1,1},{1}}
        \end{tikzpicture}};
            \node at (15,-6) {    \begin{tikzpicture}[scale=0.20] 
\planepartition{{1,1,1},{1},{1}}
        \end{tikzpicture}};
        
            \node at (0.4,-9) {    \begin{tikzpicture}[scale=0.20] 
\planepartition{{3,1},{1}}
        \end{tikzpicture}};
            \node at (3,-9) {    \begin{tikzpicture}[scale=0.20] 
\planepartition{{2,1,1},{1}}
        \end{tikzpicture}};
            \node at (6,-9) {    \begin{tikzpicture}[scale=0.20] 
\planepartition{{2,1},{1},{1}}
        \end{tikzpicture}};
            \node at (9,-9) {    \begin{tikzpicture}[scale=0.20] 
\planepartition{{2,2},{1}}
        \end{tikzpicture}};
            \node at (12,-9) {    \begin{tikzpicture}[scale=0.20] 
\planepartition{{2,1},{2}}
        \end{tikzpicture}};
            \node at (15,-9) {    \begin{tikzpicture}[scale=0.20] 
\planepartition{{2,1},{1,1}}
        \end{tikzpicture}};
        \end{tikzpicture}
        \end{tabular}
        \caption{Plane partitions of corner-hook volume and size 3,4,5.}
\label{fig:examplepim}
    \end{figure}
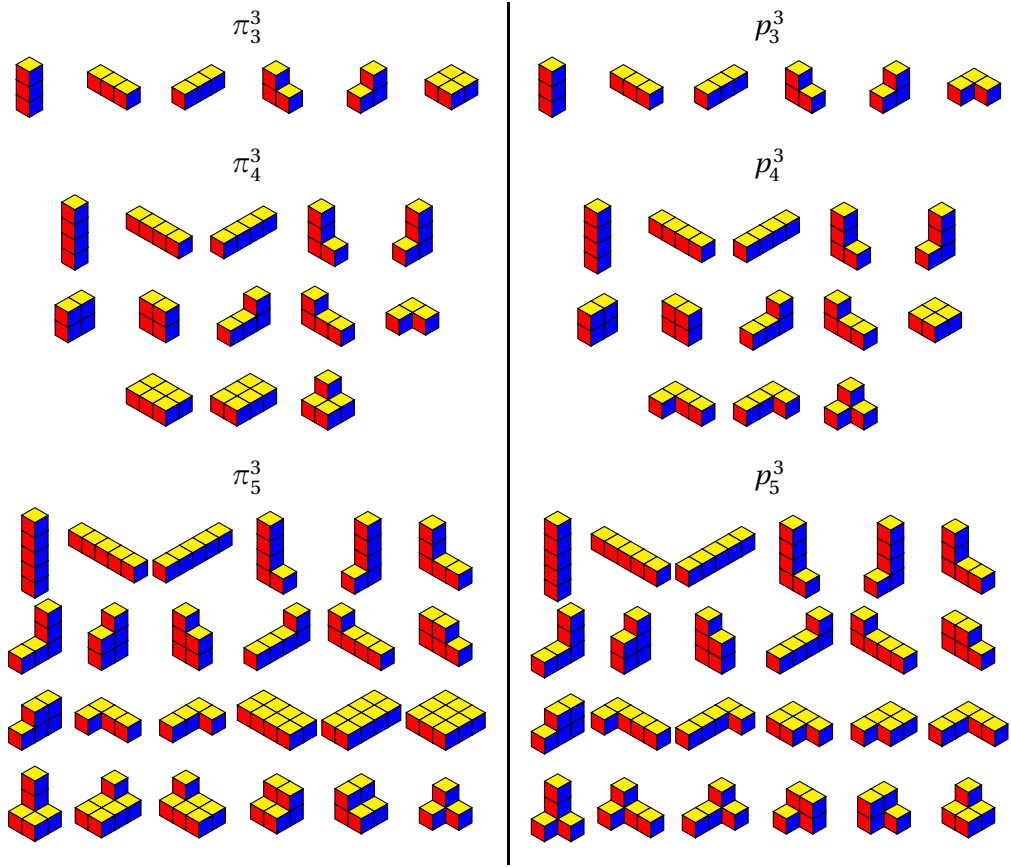
    \end{remark}

We call `MacMahon's discrepancy' the quantity
\[
e_d^k=\overline{y}_{d}^k- {y}_{d}^k,
\]
also considered in \cite{Atkin}. In \cite[Conj. 9.8]{MOTIVEMOTIVEMOTIVE} we conjectured the existence of rational generating functions for the sequences $(e_{i+k+2}^{i+1})_{i\ge 1}$, for $k\geqslant 1$. We now propose a stronger version which relates to a question raised in \cite{Atkin}, see \Cref{prop:definite-positivity}. 

\begin{conjecture}
\label{conj: Andrews}   
For every $k \in \BZ_{> 0}$ there exists a polynomial $\overline{\mathsf y}_{k}(t)\in\BQ[t]$ of degree at most $\binom{k+4}{2}- 7-k$, such that
   \[
\sum_{i\geqslant 0} \overline{y}_{i+k+2}^{i+1} t^i= \frac{ \overline{\mathsf y}_{k}(t)}{ \prod_{i=1}^k\prod_{j=1}^i (1-jt)}.
   \]  
   Moreover, for every $0\leqslant k\leqslant d-1$ we have $\overline{y}_{d}^k\geqslant 0$.
\end{conjecture}
 
\begin{remark}
\label{rem:consequences}
The nonnegativity of $\overline{y}^\bullet_\bullet$ would follow from a combinatorial interpretation of \eqref{eqn: ybarrrre}. Indeed, one could subdivide the partitions counted by $\pi_\bullet^\bullet$ according to their embedding dimension, as suggested by \Cref{fig:examplepim}. However, it is not clear whether in general this subdivision is consistent with \eqref{eqn: ybarrrre}. In particular, if such a combinatorial interpretation exists, we expect it to be possible to define integers $\overline{c}^\bullet_\bullet$ refining $\overline{y}^\bullet_\bullet$, yielding a proof of \Cref{conj: Andrews}, analogously to \Cref{thm: intro Y}.

The latter expectation comes from the fact that  the rational function  $ \prod_{j=1}^i (1-jt)^{-1}$ is the generating series of the Stirling numbers of the second kind $\left\{\begin{smallmatrix}
        n\\i
\end{smallmatrix}\right\}$, i.e.~the number of partitions of a set of $n$ elements into $i$ nonempty parts \cite[App.~A.8]{FS_book}. 
\end{remark}

\begin{remark}
\Cref{conj: Andrews} is strictly stronger than \cite[Conj. 9.8]{MOTIVEMOTIVEMOTIVE}, and is more suited for computational checks for $k,i$ arbitrarily large. Using Mathematica \cite{Mathematica} we have checked \Cref{conj: Andrews}, for $k=1,\ldots,18$ and $i<10^3$, thus confirming, thanks to  \ancillary, \cite[Conj.~9.8]{MOTIVEMOTIVEMOTIVE} as well, in the same range.
\end{remark}

Another motivation for \Cref{conj: Andrews} is a question raised in \cite{Atkin}, regarding the positivity of the numbers $e^k_d$. 
\begin{prop}
\label{prop:definite-positivity}
  If  \Cref{conj: Andrews}  holds, then $e_d^k> 0$ for $d\gg 0$. 
\end{prop}
\begin{proof}
This follows from \cite[Thm.~3.4]{GUBBIOTTI} (see also the references therein) on the asymptotic growth of the coefficients of the Maclaurin expansion of a rational generating function in terms of its radius of convergence.
\end{proof}

Via a direct check,\footnote{We stress that, thanks to \Cref{eqn: ybarrrre} and \Cref{thm: intro Y}, this check can easily be performed for $j$ arbitrarily large.} we were able to confirm the positivity of $(e^{j+1}_{j+k+2})_{j\geqslant 0}$, for $k\leqslant 18$ and $j \leqslant 10^3$. We remark that these checks are highly nontrivial, as the coefficients of $\mathsf M_{k}(t)$ (i.e.~the numerator of the generating function for $(e^{j+1}_{j+k+2})_{j\geqslant 0}$) from \cite[Conj.~9.8]{MOTIVEMOTIVEMOTIVE} are in general both positive and negative. 
 
\smallbreak
Another possible measure of the MacMahon's discrepancy is given in terms of the difference between the correct and predicted  exponents in equation \eqref{eqn: MacMahon conj intro}.
Define numbers $\omega_{m}^n,\, \overline{\omega}_m^n $ by 
\begin{align*}
  \sum_{d\geqslant 0}p_d^{n}t^d&=\prod_{m>0} \left(1-t^m\right)^{-\omega_{m}^n},\\
    \overline{\omega}_m^n &=\binom{m+n-3}{n-2},
\end{align*} 
and set
\[
\varepsilon_m^n=\overline{\omega}_m^n- {\omega}_m^n.
\]
Notice that $p_d^n$ is a polynomial of degree $d-1$ in $n$ by \Cref{prop: inversion}. As a consequence, via a direct check we see that $\omega_m^n$ is a polynomial in $n$ as well, of degree at most $m-1$. Similarly, by definition,  $\overline{\omega}_m^n $ is a polynomial in $n$ of degree $m-1$. Therefore, the `error' $\varepsilon_m^n$ is a polynomial in $n$ of degree at most $m-1$.
\begin{conjecture}[{\cite[Conj.~9.6]{MOTIVEMOTIVEMOTIVE}}]
\label{conj:erroreMac}
For every $m\geqslant 6 $ there exists an \emph{irreducible} polynomial $r_m(t)\in \BQ[t]$  of degree at most $m-6$ such that
\[
\varepsilon_m^n=\binom{n}{4}r_m(n)
\]
for all $n\geqslant 1$. In particular, the error $\varepsilon_m^n$ has degree at most $m-2$. 
\end{conjecture}
We have used \ancillary\   and Macaulay2 \cite{M2}\,  to provide checks of \cref{conj:erroreMac}.
\begin{prop}
\Cref{conj:erroreMac} holds for $m\leqslant 30$.
\end{prop}

Notice that the degree of $\varepsilon_m^n$ is predicted to be strictly smaller than $m-1$. This in particular would imply that for all $m\geqslant 1$ one has
\begin{align*}
    \lim_{n\to \infty}\frac{\overline{\omega}_m^n}{\omega_m^n}=1.
\end{align*}
\subsection{Sparsity} We conclude by formulating a conjecture on the sparsity of the numbers $p_d^n$. 

One may ask whether there are repetitions among the numbers $p_d^n$, i.e.~if there are pairs $(n,d)\not=(m,e)$ such that $p_d^n=p_{e}^m$. Since we have $1=p_1^n$ and $n=p_{2}^{n}$, the question becomes nontrivial for $d\geqslant 3$. Inspecting \ancillary, we found the following repetitions 
\begin{table}[h]
    \centering
 \begin{align*}
    p_3^5&=p_7^2=15&
    p_3^9&=p_4^5=45\\
    p_3^{14}&=p_4^7=105&
    p_3^{15}&=p_5^5=120\\
    p_3^{21}&=p_{16}^2=231&
    p_3^{65}&=p_{8}^5=2145\\
    p_3^{70}&=p_{13}^3=2485&
    p_3^{89}&=p_{4}^{27}=4005\\
    p_3^{185}&=p_{4}^{45}=17205&
    p_3^{1418}&=p_{5}^{66}=1006071\\
    p_3^{3382}&=p_{4}^{323}=5720653&
    p_3^{4898}&=p_{4}^{414}=11997651\\
    p_3^{2160324}&=p_{4}^{24100}=2333500972650 
\end{align*}  
\caption{Solutions to $p_d^n=p_{e}^m$, for $3\leqslant d <e\leqslant 30$ and  $p_e^m\leqslant 10^{18}$.}
\label{tab:repetitions}
\end{table}
which remarkably only involve, on one side of each identity, the integer $d=3$. Recall that the numbers $p_3^n$ have a peculiar combinatorial interpretation, being the number of subsets of $\Set{1,\ldots,n}$ of cardinality 1 or 2. This interpretation is not avaliable for $d\geqslant 4$. 

It is worth mentioning that, when $d=3$, for any $e>3$, the equation 
\begin{equation}
    \label{eq:hyperelliptic}
    \binom{n+1}{2}=p_e^m,
\end{equation}
defines a hyperelliptic curve with a given integral point $(1,1)$. By \cite{MR718935} the number of rational solutions of \eqref{eq:hyperelliptic} is finite, and there are algorithms  to find other integral solutions, see \cite[Sec.~10.4.1]{Cohen} and also \cite{MR2457355}. In the special case $e=4$ equation \eqref{eq:hyperelliptic} defines a smooth elliptic curve in this case there are finitely many integral points which we list in \Cref{tab:repetitions}, see \cite[Ch.~5]{SilvermanTate} for more details. It is worth mentioning that in case $e=5$ we get a singular elliptic curve. Up to our knowledge, in the case $e>d>3$ there is no general algorithm to compute all integral points.

Motivated by these observations, we propose the following conjecture.
 
\begin{conjecture}
\label{conj: spars}
Fix  $2<d<e$. If $p_d^n=p_{e}^m$ for some $n,m$, then $d=3$.
\end{conjecture}

A direct computer check using the data in \ancillary\ allowed us to confirm the conjecture in a wide range.

\begin{prop}
\Cref{conj: spars} holds in the following cases:
\begin{itemize}
    \item $e \leqslant 30$ and $m $ such that $p_e^m \leqslant  10^{23}$,
    \item $m \leqslant 3$ and $e $ such that $p_e^m \leqslant  10^{23}$.
\end{itemize}
\end{prop} 

\bibliographystyle{amsplain-nodash}
\bibliography{The_Bible}

\bigskip
\noindent
{\small{Michele Graffeo \\
\address{SISSA, Via Bonomea 265, 34136, Trieste (Italy)} \\
\href{mailto:mgraffeo@sissa.it}{\texttt{mgraffeo@sissa.it}}
}}

\bigskip
\noindent
{\small{Sergej Monavari \\
\address{\'Ecole Polytechnique F\'ed\'erale de Lausanne (EPFL),  CH-1015 Lausanne (Switzerland)} \\
\href{mailto:sergej.monavari@epfl.ch}{\texttt{sergej.monavari@epfl.ch}}
}}

\bigskip
\noindent
{\small{Riccardo Moschetti \\
\address{Università di Torino, Via Carlo Alberto 10, 10123, Torino (Italy)} \\
\href{mailto:riccardo.moschetti@unito.it}{\texttt{riccardo.moschetti@unito.it}}
}}

\bigskip
\noindent
{\small Andrea T. Ricolfi \\
\address{SISSA, Via Bonomea 265, 34136, Trieste (Italy)} \\
\href{mailto:aricolfi@sissa.it}{\texttt{aricolfi@sissa.it}}}
\end{document}